\documentclass{amsart}
\usepackage{amsmath,amsthm,amsfonts,amssymb,amscd,latexsym, mathabx}
\usepackage{graphicx,tikz}
\usepackage{pgf,pgfplots}
\usepackage{enumerate}
\usepackage{arcs}
\DeclareFontFamily{OMX}{yhex}{}
\DeclareFontShape{OMX}{yhex}{m}{n}{<->yhcmex10}{}
\DeclareSymbolFont{yhlargesymbols}{OMX}{yhex}{m}{n}
\DeclareMathAccent{\wideparen}{\mathord}{yhlargesymbols}{"F3}
\usepackage[spanish,english]{babel}
\usepackage[utf8x]{inputenc}
\usepackage[center,footnotesize]{caption2}

\newtheorem{theo}{Theorem}

\theoremstyle{definition}
\newtheorem{rk}{Remark}[section]

\def\R{\mathbb{R}}

\def\rea{\mathbb{R}}

\numberwithin{equation}{section}
\makeatletter
\newenvironment{abstracts}{%
  \ifx\maketitle\relax
    \ClassWarning{\@classname}{Abstract should precede
      \protect\maketitle\space in AMS document classes; reported}%
  \fi
  \global\setbox\abstractbox=\vtop \bgroup
    \normalfont\Small
    \list{}{\labelwidth\z@
      \leftmargin3pc \rightmargin\leftmargin
      \listparindent\normalparindent \itemindent\z@
      \parsep\z@ \@plus\p@
      
      \itemsep\medskipamount
    }%
}{%
  \endlist\egroup
  \ifx\@setabstract\relax \@setabstracta \fi
}

\newcommand{\abstractin}[1]{%
  \otherlanguage{#1}%
  \item[\hskip\labelsep\scshape\abstractname.]%
}
\makeatother

\begin{document}
\title[Coupled logistic map: A review and numerical facts]{Coupled logistic 
map: A review and numerical facts}
\author[N. Romero]{Neptal\'i Romero}
\address{Universidad Centroccidental
Lisandro Alvarado. Departamento
de Matem\'atica. Decanato de Ciencias y Tecnolog\'{\i}a.
Apartado Postal 400.
Barquisimeto, Venezuela.}
\email{nromero@ucla.edu.ve}
\author[J. Silva]{Jes\'us Silva}
\address{Universidad Centroccidental
Lisandro Alvarado. Departamento
de Matem\'atica. Decanato de Ciencias y Tecnolog\'{\i}a.
Apartado Postal 400.
Barquisimeto, Venezuela.}
\email{jesus.silva@ucla.edu.ve}

\author[R. Vivas]{Ram\'on Vivas}
\address{Universidad Nacional Experimental Polit\'ecnica Antonio José de Sucre.
Vicerrectorado de Barquisimeto.
Departamento de Estudios Básicos. Sección de Matem\'atica. 
Barquisimeto, Venezuela.} 
\email{ramon.alberto.vivas@gmail.com}


\begin{abstracts}
\abstractin{english}
This paper has a double goal, the first one is to make a slight survey of some 
theoretical results about the existence of positively invariant curves that 
allow to describe important properties of the set of bounded orbits and its 
boundary in the context of coupled logistic map. The 
second goal, in the same context, is to show a collection of figures 
obtained from computational simulation that reveal the complexity of that 
dynamical systems. 
%
%
\end{abstracts}

\maketitle

\section{Introduction}
Since a few decades ago coupled map lattices have had special interest for
physicists and scientists working in nonlinear pure and 
applied mathematics. The main reason for considering
this type of dynamical systems is that they allow to describe time
evolution of reaction-diffusion systems, population
dynamics and a huge variety of important problems in several disciplines 
as physics, chemistry, economy and even sociology. One of the most popular 
coupled map lattices is the so-called {\em coupled logistic map} (on two sides):
\begin{equation}\label{quad}
F_{\mu,\epsilon}(x,y)=((1-\epsilon)f_\mu(x)+\epsilon f_\mu(y),
(1-\epsilon)f_\mu(y)+\epsilon f_\mu(x)),
\end{equation}
where $x,y\in\R$, $f_\mu:\R\to\R$, with $f_\mu(t)=\mu t(1-t)$ is the 
logistic map and
$\epsilon$ is the coupling strength. This two-parameter family is a singular
class of coupled map lattices, it was introduced in 
independent works by Ian Fr{\o}yland \cite{froiland}, Kunihiko Kaneko 
\cite{kaneko00}, Raymond Kapral and Sergey Kuznetsov; see \cite{buni}, 
\cite{lind}. These kind of discrete dynamical systems have been intensively 
used to model a wide number of spatio-temporal phenomena in extended 
systems; see for example \cite{kaneko1}.

\smallskip

Much of the numerical and theoretical studies of \eqref{quad} obey certain 
physical interpretations under the following restrictions: 
$0\leq x,y\leq 1$, $1<\mu\leq 4$ and $0<\epsilon<1$. In this setting 
several numerical reports on the dynamics 
of \eqref{quad} are well known, cf. \cite{lind} where an important survey on 
the subject is discussed; see also \cite{D1}, \cite{D2} and 
\cite{fernandez} where some theoretical results are described. 

\smallskip

In this paper we do not consider restrictions on the state space, that 
is $x,y\in\R$; in addition we assume $\mu>1$ and $\epsilon>0$. 
Although 
it is possible that there are no physical interpretations, 
these considerations
seem interesting from the mathematical point of view. The principal aim in this 
paper is to report some computer simulations that could lead to theoretically
prove some phenomena that occur, or may be involved, in the fractalization of 
the basin of attraction of infinity, which is an attractor for the family 
$F_{\mu,\epsilon}$. Despite the fact that in the study of complicated chaotic 
behaviors, even in very simple dynamical systems, the numerical experimentation 
and their interpretations are not conclusive, they constitute an important 
support for a further analytical examination.

\smallskip

For the purpose of this paper we have taken as starting point some results 
in \cite{rsv} and \cite{rrv}, especially those 
related to the existence of invariant Jordan curves having direct
relationship with the set of points with bounded orbit; this will be 
recalled in next section. In Section 3 we show some computational simulations 
experiments followed by some interpretations. 

\section{The starting point}
We begin with some comments about the parameter space 
$(\mu,\epsilon)$. First we recall that:
$F_{\mu,\epsilon}$ and $F_{2-\mu,\epsilon}$ are topologically conjugated
when $0<\mu<1$, and the self-maps
$F_{\mu,\epsilon}$ and $F_{\mu,1-\epsilon}$ are dynamically equivalents; see 
\cite[Section 2]{rsv}. On the other hand:
$F_{0,\epsilon}$ is constant, $F_{\mu,0}$ maps $\R^2$ into the diagonal and the 
dynamics of $F_{1,\epsilon}$ is just simple: every point with bounded orbit has 
as $\omega$-limit set the origin (cf. Corollary 2.1 of \cite{rrv}). After 
these facts we fix our attention on 
$$
\epsilon\in (-\infty,0)\cup (0,1/2)\,\text{ and $\mu>1$}.
$$

We say that the coupled logistic map $F_{\mu,\epsilon}$ has {\em small 
strength} when $\epsilon\in (0,1/2)$, and it has {\em large strength} if 
$\epsilon\in (-\infty,0)$.
Regardless of strength type, there are some common dynamical properties; we 
refer to \cite{rrv} and \cite{rsv} for their proofs:

\smallskip
\noindent
$\bullet$
There are always two fixed points: $O=(0,0)$ and 
$P_\mu=\left(\frac{\mu-1}{\mu}, \frac{\mu-1}{\mu}\right)$, the nature of these 
points depends on $\mu$ and $\epsilon$. Two other fixed points appear 
in a certain region 
of the parameter space $(\mu,\epsilon)$, although these points have the same 
algebraic expression for the different strengths:
$R(p_{\mu,\epsilon})$ and $p_{\mu,\epsilon}=(p_-,p_+)$, where $R$ is the 
reflection $R(x,y)=(y,x)$ and
$$
p_\pm=\frac{k\mu\pm\sqrt{2(\mu-1)\mu k-\mu^2 k^2}}{2\mu}\text{ and }
k=1-\frac{1}{\mu(1-2\epsilon)},
$$
they appear through a pitchfork bifurcation for the fixed points $O$ (small 
strength) and $P_\mu$ (large strength).

\smallskip
\noindent
$\bullet$
The infinity is always an attractor; that is, there exists a compact 
neighborhood whose complement $U$ is positively invariant under 
$F_{\mu,\epsilon}$ (i.e. $F_{\mu,\epsilon}(U^c)\subset U^c$) and
$\|F_{\mu,\epsilon}^k(z)\|\to +\infty$ as $k\to +\infty$ for all $z\in U^c$. 
It is easy to see that if $C$ is the circle given by $x^2+y^2=x+y$, 
then the closed disk with boundary $C$ takes the role of the 
set $U$ above. Thus, if $ext\,C$ denotes the unbounded component of 
the complement of $C$, then it is contained in
$$
B_\infty(F_{\mu,\epsilon})=\{z\in\R^2:\lim_{k\to +\infty}\|F_{\mu,\epsilon}
^k(z)\|= +\infty\},
$$ 
which is called {\em basin of attraction of $\infty$}. It is also easy to see 
that
$$
B_\infty(F_{\mu,\epsilon})=\bigcup_{n\geq 0}F_{\mu,\epsilon}^{-n}(ext\,C)\,
\text{ and }\,
B^c_\infty(F_{\mu,\epsilon})=\bigcap_{n\geq 0}F_{\mu,\epsilon}^{-n}(cl(int\,C)),
$$
where $int\,C$ is the bounded component of the complement of $C$, $cl(A)$ and 
$A^c$ denote, respectively, the closure and complement of the set $A$. 
Obviously $B^c_\infty(F_{\mu,\epsilon})$ is the set of points with bounded 
orbit. The two identities above are satisfied 
when $C$ is substituted by any Jordan curve (i.e. a simple and closed curve) 
$\Gamma$ such that 
$ext\,\Gamma\subset B_\infty(F_{\mu,\epsilon})$; here $ext\,\Gamma$ has the 
same meaning of $ext\,C$.

\smallskip
\noindent
$\bullet$
The diagonal $\Delta=\{(x,x):x\in\R\}$ is always an invariant set for 
$F_{\mu,\epsilon}$ and the dynamics on this set is the same of the logistic 
map $f_\mu$. In particular, when $\mu>4$ the restriction 
$F_{\mu,\epsilon}\vert_{\Delta}$ of $F_{\mu,\epsilon}$ to $\Delta$ has an 
invariant Cantor set $K_\mu$ as its nonwandering set and its dynamics is 
topologically conjugated to the 2-symbols unilateral shift. On the other hand, 
when $\mu>4$ is large enough, the set of 
points with bounded orbits $B^c_\infty(F_{\mu,\epsilon})$ is also a Cantor set, 
containing properly $K_\mu$, and the dynamics of $F_{\mu,\epsilon}$ on 
that set is topologically conjugated to the unilateral shift on four symbols.

\subsection{Critical set and image set}
In the study of dynamical systems provided 
by differentiable and non-invertible transformations the set of critical 
points and the set of critical values have a relevant role. We recall that for a 
differentiable 
endomorphisms the set of critical points 
is the set of points where the Jacobian matrix is singular. For 
$F_{\mu,\epsilon}$ its set of critical points is given by
$$
\ell=\ell_1\cup\ell_2:=
\{(\tfrac{1}{2},y):y\in\rea\}\cup \{(x,\tfrac{1}{2}):x\in\rea\};
$$
and its set of critical values is 
$L:=F_{\mu,\epsilon}(\ell)=L_1\cup L_2$, where
\begin{gather*}
L_1=F_{\mu,\epsilon}(\ell_1)=
\{(x,y): y=\tfrac{1-\epsilon}{\epsilon}x-\tfrac{(1-2\epsilon)\mu}{4\epsilon},\,
x\leq \tfrac{\mu}{4}\}\,\text{ and}\\
L_2=F_{\mu,\epsilon}(\ell_2)=
\{(x,y): 
y=\tfrac{\epsilon}{1-\epsilon}x+\tfrac{(1-2\epsilon)\mu}{4(1-\epsilon)},\,
x\leq \tfrac{\mu}{4}\}.
\end{gather*}
\begin{figure}[h!]
\centering
\includegraphics[scale=.85]{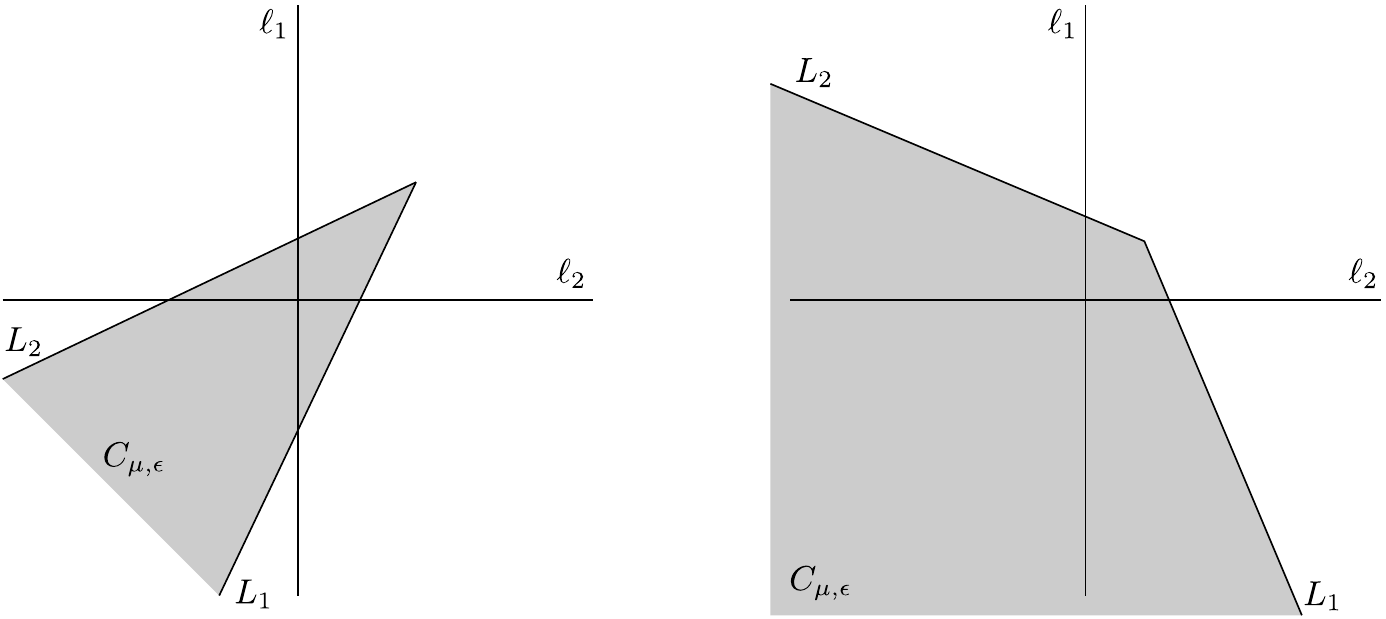}
\caption{Illustration of the critical set and critical values of 
$F_{\mu,\epsilon}$. Left figure corresponds to small strength, the other 
one to the large strength; regions in shadow represent 
$C_{\mu,\epsilon}$.}
\end{figure}
By using these rays is defined the cone 
$$
C_{\mu,\epsilon}=\left\{(x,y):
\tfrac{1-\epsilon}{\epsilon}x-\tfrac{(1-2\epsilon)\mu}{4\epsilon}
\leq y \leq
\tfrac{\epsilon}{1-\epsilon}x+\tfrac{(1-2\epsilon)\mu}{4(1-\epsilon)},\,
x\leq\tfrac{\mu}{4}\right\}.
$$
It is easy to check that
$F_{\mu,\epsilon}(\R^2)=C_{\mu,\epsilon}$;
moreover,
every point in the interior of cone $C_{\mu,\epsilon}$ has four preimages,
they are located
symmetrically respect to $\ell_1$ and $\ell_2$; indeed,
$F_{\mu,\epsilon}^{-1}(z,w)=\{(x_\pm(z,w),y_\pm(z,w))\}$, where
\begin{equation}\label{pre}
\begin{gathered}
x_\pm(z,w)=\tfrac{1}{2}\left[1\pm\sqrt{1+\tfrac{4(\epsilon
w-(1-\epsilon)z)}{\mu(1-2\epsilon)}} \right]\,\text{ and }\\ 
y_\pm(z,w)=\tfrac{1}{2}\left[1\pm \sqrt{1+\tfrac{4(\epsilon z
-(1-\epsilon)w)}{\mu(1-2\epsilon)}}\right];
\end{gathered}
\end{equation}
points outside $C_{\mu,\epsilon}$ have no preimages, and
the restriction of $F_{\mu,\epsilon}$ to $\ell_i$ is two-to-one onto $L_i$, 
$i=1,2$. In particular, if $\gamma$ is an injective curve joining $L_1$ and 
$L_2$ with extreme points in these rays and the others ones are in the 
interior of $C_{\mu,\epsilon}$, then $F_{\mu,\epsilon}^{-1}(\gamma)$ is a 
Jordan curve symmetric with respect to the critical lines $\ell_1$ and $\ell_2$ 
and surrounding the critical point $c=(1/2,1/2)$.

\begin{figure}[h!]
\centering
\includegraphics[scale=1]{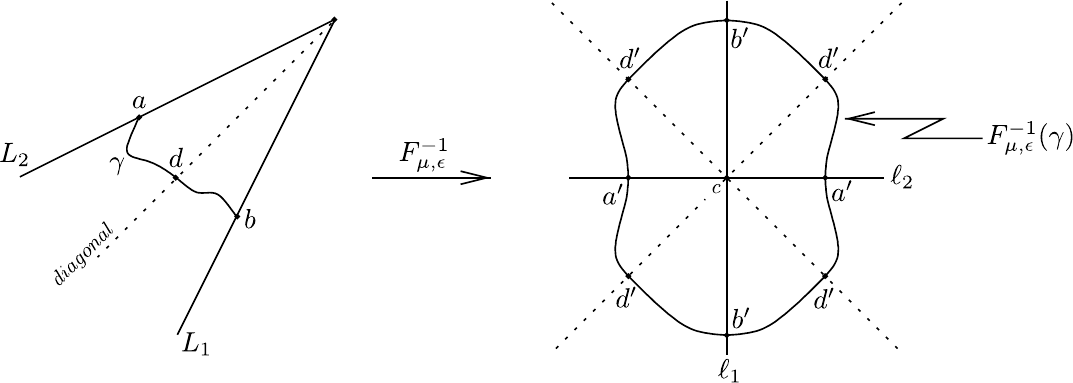}
\caption{Graphical illustration of how the preimage of a curve $\gamma$ inside 
$C_{\mu,\epsilon}$, with extremes at $L_1$ and $L_2$, is arranged. 
The points 
$a'$, $b'$ and $d'$ are, respectively, preimages of $a,b$ and $d$.}
\end{figure}

\subsection{Small strength: $\boldsymbol{\epsilon\in (0,1/2)}$}
\label{small1}
In this case there are three curves playing an important meaning in the 
dynamic description of $F_{\mu,\epsilon}$, such a curves 
are:
$$
\epsilon \longmapsto \mu_0(\epsilon)=\tfrac{1}{1-2\epsilon},\,
\epsilon \longmapsto \mu_1(\epsilon)=\tfrac{4(1-\epsilon)}{1-2\epsilon}
\,\text{ and }
\epsilon \longmapsto 
\mu'(\epsilon)=1+\sqrt{\tfrac{3-2\epsilon}{1-2\epsilon}};
$$
the functions $\mu_0$ and $\mu_1$ have the interval $(0,1/2)$ as domain, while 
the 
domain of $\mu'$ is $(0,3/8]$.
The first curve ($\epsilon\longmapsto\mu_0(\epsilon)$) defines the locus where 
the 
change in the hyperbolic nature of 
the fixed point at the origin is marked: $O$ is a hyperbolic saddle when 
$\mu<\mu_0(\epsilon)$ and it 
is a repeller if $\mu>\mu_0(\epsilon)$. At $\mu=\mu_0(\epsilon)$ the origin has 
a pitchfork bifurcation and the fixed points $R(p_{\mu,\epsilon})$ and 
$p_{\mu,\epsilon}$ are correctly defined if and only if 
$\mu>\mu_0(\epsilon)$. 
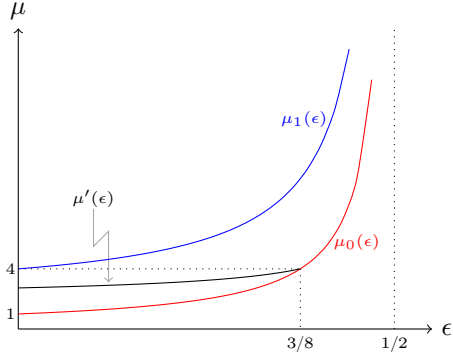
\begin{figure}[h!]
\centering
\begin{tikzpicture}[smooth, xscale=10, yscale=.2]
\draw[->] (0,0) -- (.55,0) node[right] {$\epsilon$};
\draw[->] (0,0) -- (0,20) node[above] {$\mu$};
\draw[dotted] (0.375,4) -- (0.375,0);
\draw[dotted] (0,4) -- (0.375,4);
\draw[->,thin,color=gray!70!] (0.1,8) -- (0.1,5.5) -- (0.12,6.5) -- (0.12,3.1); 
\draw[color=red]
plot[domain=0:.47] (\x,{1/(1-2*\x)});
\draw[color=blue]
plot[domain=0:.44] (\x,{4*(1-\x)/(1-2*\x)});
\draw
plot[domain=0:.375] (\x,{1+(sqrt{(3-2*\x)}/sqrt{(1-2*\x)})});
\draw[dotted] (.5,0)--(.5,20);
\node at (.5,-1) {\tiny{$1/2$}};
\node at (.375,-1) {\tiny{$3/8$}};
\node at (-.01,1) {\tiny{$1$}};
\node at (-.01,4) {\tiny{$4$}};
\node[color=red] at (.45,5.8) {\tiny{$\mu_0(\epsilon)$}};
\node[color=blue] at (.38,14) {\tiny{$\mu_1(\epsilon)$}};
\node at (.1,8.7) {\tiny{$\mu'(\epsilon)$}};
\end{tikzpicture}
\caption{Graphical illustration of the functions $\epsilon \longmapsto 
\mu_0(\epsilon)$, $\epsilon \longmapsto \mu_1(\epsilon)$ and 
$\epsilon\longmapsto \mu'(\epsilon)$.}
\label{3curvas1}
\end{figure}

From the second curve ($\epsilon\longmapsto\mu_1(\epsilon)$) is ensured the 
existence of a Jordan curve through which 
is possible to make a description of the basin of attraction of $\infty$ and 
its boundary $\partial B_\infty(F_{\mu,\epsilon})$. Through this curve is also 
characterized the location of the intersection point $(q_{\mu,\epsilon},0)$ 
between $L_1$ and the horizontal axis, here 
$q_{\mu,\epsilon}=\frac{(1-2\epsilon)\mu}{4(1-\epsilon)}$. Indeed, 
$q_{\mu,\epsilon}>1$ if and only if $\mu>\mu_1(\epsilon)$.

\smallskip

The third curve ($\epsilon\longmapsto\mu'(\epsilon)$) is related 
to the {\em synchronization} of orbits, an important 
dynamic phenomenon that has captured the attention of several authors. 
In addition, the first curve is also
a reference to the description of the synchronized points. 
We denote by 
$\mathcal{S}(\mu,\epsilon)$ the set of {\em synchronized points} of
$F_{\mu,\epsilon}$; that is, the set of points $(x,y)\in\rea^2$
whose positive orbit is bounded and 
$\lim_{n\to+\infty}|x_n-y_n|=0$, where $(x_n,y_n)=F_{\mu,\epsilon}^n(x,y)$
for all $n\geq 1$.

\smallskip

Now we summarize the principal results in \cite{rrv}.

\begin{theo}\label{a}
If $\mu>1$ and $\epsilon\in (0,1/2)$ satisfy $\mu\leq\mu_1(\epsilon)$, then 
there 
exists a positively invariant Jordan curve $\Gamma$ containing 
$F_{\mu,\epsilon}^{-1}(O)$ such that:

\begin{enumerate}[a)]
\item 
The open set $ext\,\Gamma$ is contained in $B_\infty(F_{\mu,\epsilon})$ and 
$\Gamma=\partial B^o_\infty(F_{\mu,\epsilon})$; in particular
$$
B_\infty(F_{\mu,\epsilon})=\bigcup_{n\geq 0}F_{\mu,\epsilon}^{-n}(ext\,\Gamma)\,
\text{ and }\,
B^c_\infty(F_{\mu,\epsilon})=\bigcap_{n\geq 
0}F_{\mu,\epsilon}^{-n}(cl(int\,\Gamma)).
$$
The symbol $B^o_\infty(F_{\mu,\epsilon})$ 
means the unbounded connected component of $B_\infty(F_{\mu,\epsilon})$ 
(i.e. the immediate basin of $\infty$) and $\partial 
B^o_\infty(F_{\mu,\epsilon})$ denotes its boundary.

\item
If $\mu\leq 4$, then $B_\infty(F_{\mu,\epsilon})$ is connected; in 
this case $B_\infty(F_{\mu,\epsilon})=ext\,\Gamma$. On the other hand, when 
$4<\mu\leq\mu_1(\epsilon)$ the basin $B_\infty(F_{\mu,\epsilon})$ has 
infinitely many connected components; in this case 
$\partial B_\infty(F_{\mu,\epsilon})=
\bigcup_{n\geq 0}F^{-n}_{\mu,\epsilon}(\Gamma)$.

\item\label{a-c}
When $\mu\leq\mu_0(\epsilon)$, 
$\mathcal{S}(\mu,\epsilon)=B^c_\infty(F_{\mu,\epsilon})$. Consequently the 
$\omega$-limit set of every point in $int\,\Gamma$ is contained in $\Delta$; 
further, if $4<\mu<\mu_0(\epsilon)$, then 
$\mathcal{S}(\mu,\epsilon)$ is the union of the stable manifolds of the points 
in the Cantor set on $\Delta$. In addition, if $\mu_0(\epsilon)<\mu\leq 
\mu'(\epsilon)$, then 
$\mathcal{S}(\mu,\epsilon)=\left (B^c_\infty(F_{\mu,\epsilon})\setminus 
\Gamma\right)\cup F_{\mu,\epsilon}^{-1}(O)$.

\item
When $\mu<\mu_0(\epsilon)$, $\Gamma$ is the connected component of the stable 
manifold $W^s(O)$ of the origin containing it; moreover, 
$W^s(O)=\bigcup_{n\geq 0}F^{-n}_{\mu,\epsilon}(\Gamma)$. If 
$\mu=\mu_0(\epsilon)$, then $\Gamma$ contains the unique local center-stable 
manifold of $O$. In particular, for $\mu\leq \mu_0(\epsilon)$ the curve 
$\Gamma$ is $C^\infty$ and the $\omega$-limit set of every $z\in\Gamma$ is the 
origin.
\end{enumerate}
\end{theo}

Some comments are natural from the statements in the previous theorem.

\begin{rk}\label{rk1}

\smallskip
\noindent
1)
For the parametric region $\mu_0(\epsilon)<\mu\leq\mu_1(\epsilon)$, the curve 
$\Gamma$ is obtained as a collage of four others, one of them, denoted by 
$\Gamma_b$, has extreme points at $O$ and $S=(1,0)$, it is the graph of some 
function $\gamma:[0,1]\to\R$ which is contained in  
$D=cl(int\,C)\cap\{(x,y):y\leq 0\}$. Indeed, $\gamma$ is the fixed point of a 
contracting operator acting on the complete metric space of functions 
defined in $[0,1]$ such that each one of them has Lipschitz constant 
less or equal than 1, it is symmetric respect to $1/2$, fixes $t_0=0$, its 
graph is in $D$ and it is tangent to the antidiagonal at $t_0$. The distance in 
this space is the usual between bounded functions.
Thus, the curve $\Gamma$ is given by
$\Gamma=\Gamma_b\cup\Gamma_t\cup\Gamma_\ell\cup\Gamma_r$ where
\begin{equation}\label{Gamma}
\begin{gathered}
\Gamma_t=\{(t,1-\gamma(t)): t\in [0,1]\},\,
\Gamma_\ell=\{(\gamma(t),t): t\in [0,1]\}\\
\text{and $\Gamma_r=\{(1-\gamma(t),t): t\in [0,1]\}$};
\end{gathered}
\end{equation}
in addition, $F_{\mu,\epsilon}^{-1}(\Gamma_b)=\Gamma_b\cup \Gamma_t$,
$F_{\mu,\epsilon}^{-1}(\Gamma_\ell)=\Gamma_\ell\cup \Gamma_r$ and
$F_{\mu,\epsilon}(\Gamma)\subset \Gamma_b\cup \Gamma_\ell$; see Figure 
\ref{gamma}.
\begin{figure}[h!]
\centering
\includegraphics[scale=1]{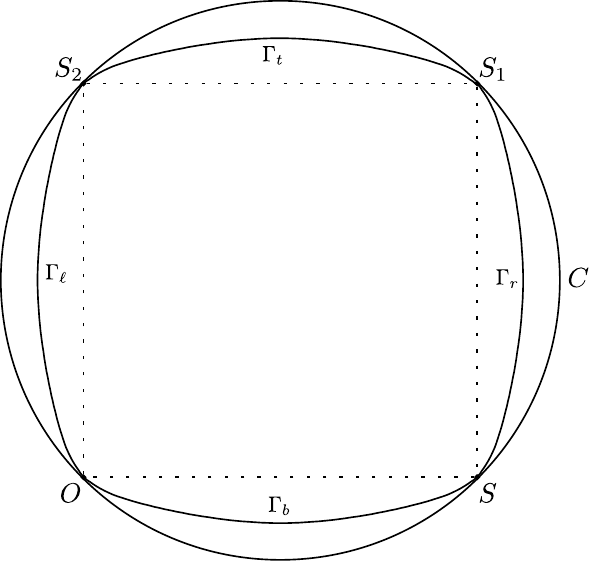}
\caption{Idealization of the Jordan curve 
$\Gamma=\Gamma_b\cup\Gamma_t\cup\Gamma_\ell\cup\Gamma_r$, which is also 
obtained as the limit curves given by the preimages of the circle $C$.}
\label{gamma}
\end{figure}
That operator is defined as follows: given a function in that space, 
the preimage under $F_{\mu,\epsilon}$ of its graph is 
union of two graphs: one of them is located above $\ell_2$, it joins the
$S_1=(1,1)$ and $S_2=(0,1)$; the other one is below $\ell_2$ and connects the 
points $O$ with $S$. Just the function defining this last graph is the image 
of the operator of the given function.

\smallskip

If $\mu$ and $\epsilon$ satisfy $1<\mu\leq\mu_0(\epsilon)$, the 
operator above is not necessarily a contraction. However, to obtain $\Gamma$ 
one takes the curve given by the arc of the circle $C$ in 
$D$, then the forward iterations of this arc under the operator above define
an increasing sequence of functions whose limit function and the preceding 
collage determine a positively invariant Jordan curve whose unbounded 
component of its complement is contained in $B_\infty(F_{\mu,\epsilon})$; the 
uniqueness of such a curve is consequence of the synchronization property due 
to the coincidence of $\Gamma$ with the stable ($1<\mu<\mu_0(\epsilon)$) and 
center-stable ($\mu=\mu_0(\epsilon)$) manifold of the origin. 
In this way, for all $\mu$ and $\epsilon$ with $1<\mu\leq\mu_1(\epsilon)$, 
the curve $\Gamma$ is the unique positively invariant Jordan curve such that 
$F^{-1}_{\mu,\epsilon}(O)\subset \Gamma$ and  
$ext\,\Gamma=B^o_\infty(F_{\mu,\epsilon})$.

\smallskip
\noindent
2)
For all $\mu$ and $\epsilon\in (0,1/2)$ with $4<\mu\leq \mu_1(\epsilon)$ 
the attraction basin $B_\infty(F_{\mu,\epsilon})$ has infinitely many connected 
components inside $int\,\Gamma$.
For this sector in the parameter space it is proved that there are two points 
$p_i\in L_i\cap \Gamma$ ($i=1,2$) such that the arc of $\Gamma$ containing 
$p_1,p_2$ and $S_1$ and the pieces of those rays between $p_1,p_2$ and 
$d_\mu=(\mu/4,\mu/4)$ determine a triangular region in $ext\,\Gamma$ in such a 
way that its preimage under $F_{\mu,\epsilon}$ defines a connected component of 
$B_\infty(F_{\mu,\epsilon})$ inside $int\,\Gamma$ whose boundary is part of 
$F^{-1}_{\mu,\epsilon}(\Gamma)$ and surrounds the critical point $c$; 
see Figure \ref{infa}. 
\begin{figure}[h!]
\centering
\includegraphics[scale=1]{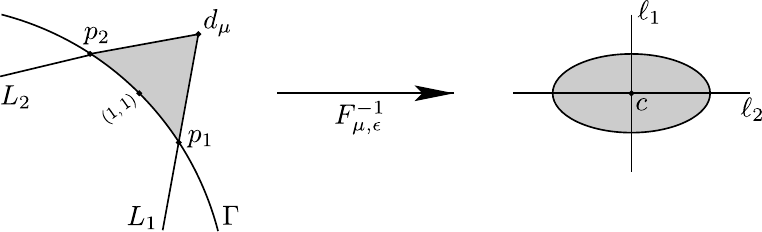}
\caption{The connect component of $B_\infty(F_{\mu,\epsilon})$ inside 
$int\,\Gamma$ (shaded region on the right side) is obtained as preimage of 
the triangular region described above.}
\label{infa}
\end{figure}
Just the recursive preimages of this component produce 
infinitely many other connected components of $B_\infty(F_{\mu,\epsilon})$ 
inside $int\,\Gamma$. Every component in the segment $OS_1$ of the
complement of the Cantor set $K_\mu$ is contained in some of those preimages 
and each of which contains at most two components of $K^c_\mu$ in that segment. 
\begin{figure}[h!]
\centering
\includegraphics[scale=1]{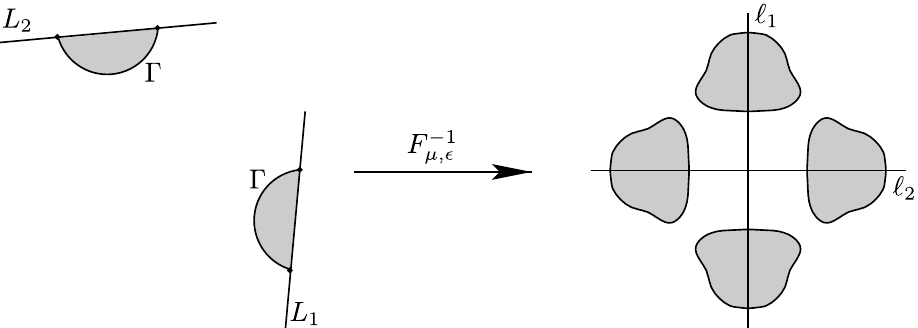}
\caption{In the left side it is illustrates a pair of sectors (shadowed 
regions) whose preimages, except their boundaries, are components of 
$B_\infty(F_{\mu,\epsilon})$ inside $int\,\Gamma$; observe that they do not 
intersect the segment $OS_1$.}
\label{infb}
\end{figure}

It is also possible that there are other components of 
$B_\infty(F_{\mu,\epsilon})$ inside $int\,\Gamma$ such that none of its images 
or preimages intersects the complement of $K_\mu$ in $\Delta$.
These other components of $B_\infty(F_{\mu,\epsilon})$ are produced only as 
preimages of sectors bounded by an arc in $\Gamma$ and a segment in 
$L_1$ (resp. $L_2$) with the same extreme points, such that each one of these 
sectors, except that arc, is contained in $ext\,\Gamma$. In view of the nice 
symmetric properties of $\Gamma$ respect to $\ell_1$ and $\ell_2$, this kind of 
sectors appear in pairs, one is the reflection of the other. The preimage of 
each one of these sectors is the disjoint union of two topological disks
displayed symmetrically and whose boundaries are contained in 
$F_{\mu,\epsilon}^{-1}(\Gamma)$; see 
Figure \ref{infb}.

\smallskip

Now consider $\mu>1$ and $\epsilon\in (0,1/2)$ such that $\mu>\mu_1(\epsilon)$, 
that is $q_{\mu,\epsilon}>1$; recall the meaning of the curve 
$\epsilon\longmapsto\mu_1(\epsilon)$ above. In this case two mutually exclusive 
situations occur: either there exists an integer $n\geq 0$ such that 
$F^{-n}_{\mu,\epsilon}(C)$ is contained into the interior of 
$C_{\mu,\epsilon}$, or for all $n\geq 0$ the connect component of 
$F^{-n}_{\mu,\epsilon}(C)$ containing $O$ also contains 
$F^{-1}_{\mu,\epsilon}(O)$ and intersect the rays $L_1$ and $L_2$. In the first 
configuration the set of points with bounded orbit has infinitely components 
and none of them contains $F^{-1}_{\mu,\epsilon}(O)$, hence it is impossible 
the existence of a Jordan as above; observe that if $\mu$ is large enough, then 
this arrangement is achieved. In the second situation, we do not know 
a theoretical result guaranteeing the existence of such a curve; however, when 
$\mu$ is close enough to $\mu_1(\epsilon)$ there are computacional simulations 
that suggest such an existence; see Figure \ref{beyond}. Observe that for 
any curve in $D$ joining $O$ and $S$ its preimage under $F_{\mu,\epsilon}$ 
is disjoint union of four arcs: two above $\ell_2$ and two below it; each one 
of these arcs connects a preimage of $O$ with the closest preimage of $S$; see 
Figure \ref{infc}.
\begin{figure}[h!]
\centering
\includegraphics[scale=1]{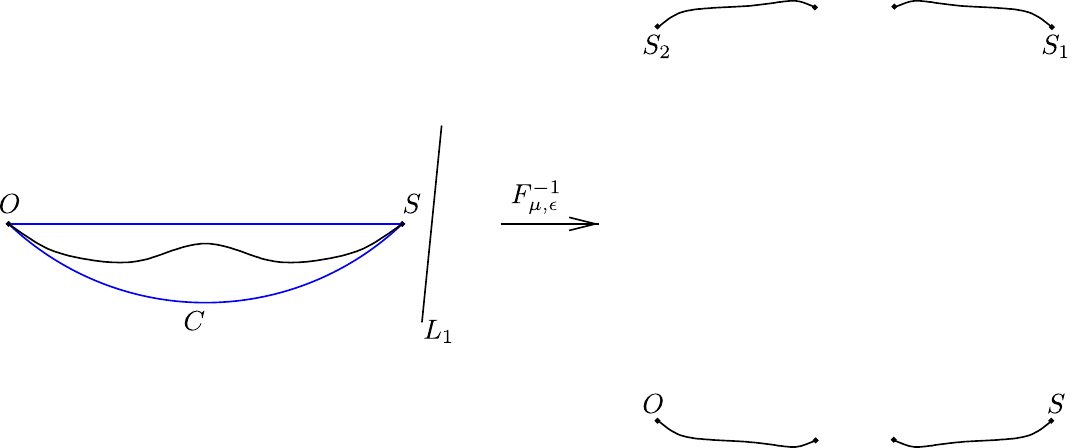}
\caption{In the right side of the figure is illustrated the preimage of the 
curve contained in $D$ (left side) when $\mu>\mu_1(\epsilon)$, i.e. when 
$q_{\mu,\epsilon}>1$.}
\label{infc}
\end{figure}
Despite this difficulty one can define a similar operator to the previous and 
then try to get, in some sense, a limit curve for the iterations of that 
operator. The ideas of such a procedure are the same as those discussed in item 
3 of Remark \ref{rk2}.
This is still under development, we hope to present 
results in a forthcoming article.
\end{rk}
\subsection{Large strength: $\boldsymbol{\epsilon<0}$}

When the parameters $\mu$ and $\epsilon$ are such that $\mu>1$ and $\epsilon<0$ 
(or equivalently $\epsilon>1$) the fixed points of $F_{\mu,\epsilon}$ have 
different nature than in the small strength case. In particular, for the 
fixed points on the diagonal, $O$ and $P_\mu$, are well known the following 
facts:
\begin{enumerate}[(i)]
\item
The origin $O$ is always a hyperbolic repeller.
\item
If $1<\mu<\mu_0'(\epsilon):=\frac{1-4\epsilon}{1-2\epsilon}$, then 
$P_\mu$ is a hyperbolic saddle.
\item
When 
$\mu_0'(\epsilon)<\mu<\mu_2(\epsilon):=
\frac{3-4\epsilon}{1-2\epsilon}$, $P_\mu$ is a hyperbolic 
attractor.
\item
If $\mu_2(\epsilon)<\mu<3$, then $P_\mu$ is again a hyperbolic saddle.
\item
For $\mu>3$, $P_\mu$ is always a hyperbolic repeller.
\end{enumerate}
\begin{center}
\begin{figure}[h!]
\centering
\begin{tikzpicture}[smooth, xscale=20, yscale=.5]
\draw[->] (0,0) -- (-.55,0) node[right] {\hspace{-.4cm}$\epsilon$};
\draw[->] (0,0) -- (0,6) node[above] {$\mu$};
\draw[dotted] (0,2) -- (-.55,2);
\draw[color=red] plot[domain=-.55:0] (\x,{(1-4*\x)/(1-2*\x)});
\draw[color=blue] plot[domain=-.55:0] (\x,{4*(1-\x)/(1-2*\x)});
\draw plot[domain=-.55:0] (\x,{(3-4*\x)/(1-2*\x)});
\node at (.01,1) {\tiny{$1$}};
\node at (.01,2) {\tiny{$2$}};
\node at (.01,3) {\tiny{$3$}};
\node at (.01,4) {\tiny{$4$}};
\node[color=blue] at (-.5,3.5) {\tiny{$\mu_1(\epsilon)$}};
\node[color=red] at (-.5,1) {\tiny{$\mu_0'(\epsilon)$}};
\node at (-.06,2.5) {\tiny{$\mu_2(\epsilon)$}};
\end{tikzpicture}
\caption{Graphs of the functions $\epsilon \longmapsto \mu_0'(\epsilon)$,
$\epsilon \longmapsto \mu_1(\epsilon)$ and 
$\epsilon \longmapsto \mu_2(\epsilon)$.}
\end{figure}
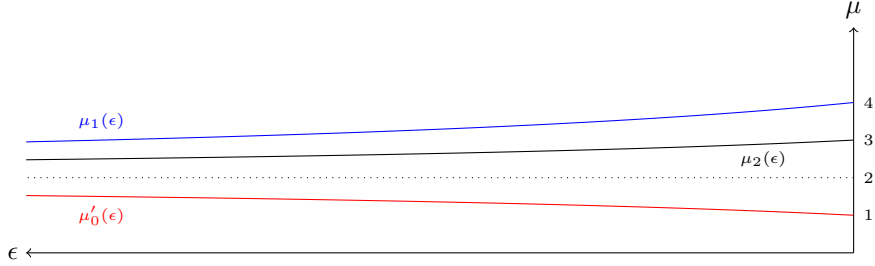
\end{center}

The fixed points outside the diagonal, $p_{\mu,\epsilon}$ and 
$R(p_{\mu,\epsilon})$, appear when $\mu>\mu_0'(\epsilon)$. 
The curve $\{(\epsilon,\mu_0'(\epsilon)): 
\epsilon<0\}$ is the locus where a pitchfork bifurcation occurs for the fixed 
point $P_\mu$. Another curve in the parameter space with a prominent 
performance is given by 
$\mu_1(\epsilon):=\frac{4(1-\epsilon)}{1-2\epsilon}$. 
As in the case of small strength, the curve 
$\{(\epsilon,\mu_1(\epsilon)): \epsilon<0\}$ determines the region in the 
parameter space where:

\begin{enumerate}[a)]
\item
The intersection point $(q_{\mu,\epsilon},0)$ between 
$L_1$ and the $x$-axis is outside the square $Q=[0,1]\times [0,1]$; indeed, 
$q_{\mu,\epsilon}\leq 1$ if, and only if, $\mu\leq \mu_1(\epsilon)$.
\item
It is guaranteed the existence of an invariant curve providing 
information about $B_\infty(F_{\mu,\epsilon})$ and its boundary, see statement 
below.
\end{enumerate}

The square $Q$ refines the first estimate of the basin of attraction of 
$\infty$:
the boundary $\partial Q$ of $Q$ is such that 
$ext\,\partial Q$ is contained in $B_\infty(F_{\mu,\epsilon})$, hence
$$
B_\infty(F_{\mu,\epsilon})=\bigcup_{n\geq 
0}F_{\mu,\epsilon}^{-n}(ext\,\partial Q)\,
\text{ and }\,
B^c_\infty(F_{\mu,\epsilon})=\bigcap_{n\geq 
0}F_{\mu,\epsilon}^{-n}(Q);
$$
in addition, $\partial Q\setminus F^{-1}_{\mu,\epsilon}(O)\subset 
B_\infty(F_{\mu,\epsilon})$.

\smallskip

Now we compile the statements of the more relevant results in \cite{rsv}.

\begin{theo}\label{b}
If $\mu>1$ and $\epsilon<0$ satisfy $\mu\leq\mu_1(\epsilon)$, then 
there 
exists a Lipschitzian and positively invariant curve $\Gamma$ containing 
$F_{\mu,\epsilon}^{-1}(O)$ and the fixed points outside the diagonal, such that:
\begin{enumerate}[a)]
\item
If $\mu_0'(\epsilon)<\mu\leq \mu_1(\epsilon)$, then $\Gamma$ is a Jordan curve 
and $ext\,\Gamma=B_\infty^o(F_{\mu,\epsilon})$; in other words, 
$ext\,\Gamma\subset B_\infty(F_{\mu,\epsilon})$ and $\Gamma=\partial 
B^o_\infty(F_{\mu,\epsilon})$. As in the small strength case
$$
B_\infty(F_{\mu,\epsilon})=\bigcup_{n\geq 0}F_{\mu,\epsilon}^{-n}(ext\,\Gamma)\,
\text{ and }\,
B^c_\infty(F_{\mu,\epsilon})=\bigcap_{n\geq 
0}F_{\mu,\epsilon}^{-n}(cl(int\,\Gamma)).
$$

\item
When $\mu\leq\mu_0'(\epsilon)$ the curve $\Gamma$ is the union of the 
straight segments $O S_1$ and $SS_2$. In this case, 
$B_\infty(F_{\mu,\epsilon})=\R^2\setminus\Gamma$ and the $\omega$-limit set of 
every point $z$ in $\Gamma\setminus F_{\mu,\epsilon}^{-1}(O)$ is $P_\mu$; that 
is $F_{\mu,\epsilon}^n(z)\to P_\mu$ if $n\to +\infty$.
\end{enumerate}
\end{theo}

\begin{rk}\label{rk2}
1) Let $T$ be the closed triangular region defined by the vertices $O,S$ and 
$c=(1/2,1/2)$. Denote by $Lip_1$ the complete metric space of all the 
Lipschitzian 
functions $h:[0,1]\to\R$ with Lipschitz constant $L(h)\leq 1$ such that: 
$h(t)=h(1-t)$ for all $t\in [0,1/2]$, $h(0)=0$ and whose graph is contained in 
$T$; the distance here the usual one between bounded functions. Given any $h\in 
Lip_1$, its graph intersects $L_1$ in only one point $(t^*,h(t^*))$ where 
$0<t^*\leq 1$; so, by using \eqref{pre}
$$
\{(x_-(t,h(t)), y_-(t,h(t))):t\in [0,t^*]\}\cup
\{(x_+(t,h(t)), y_-(t,h(t))):t\in [0,t^*]\}
$$
is a graph of some function in $Lip_1$. This fact describes an assignment 
defining an operator $\Psi$ on $Lip_1$, which  
is not necessarily a contraction; so to obtain $\Gamma$ one 
begins with the null function $h_0$ and consider its forward orbit 
$(h_n)_{n\geq 0}$, with $h_n=\Psi^n(h_0)$. Then it is proved 
that 
$h_n(t)<h_{n+1}(t)$ for all $n\geq 0$ and $t\in (0,1)$. Hence a limit function 
$\widehat{h}\in Lip_1$ is obtained, which is clearly a fixed point of $\Psi$. 
Through the expanding property of the origin it is showed that the 
graph of $\widehat{h}$, in what follows $\Gamma_b$, is tangent to the diagonal 
and the antidiagonal at $O$ and $S$, respectively; with this graph is 
constructed 
the juxtaposition $\Gamma=\Gamma_b\cup\Gamma_t\cup\Gamma_\ell\cup\Gamma_r$, 
where $\Gamma_t,\Gamma_\ell$ and $\Gamma_r$ are as in the small strength case. 
At this point, if 
$\mu_0'(\epsilon)<\mu\leq\mu_1(\epsilon)$, then is proved that 
$\widehat{h}(t)<t$ for all $t\in (0,1/2]$ and the fixed points 
$p_{\mu,\epsilon}$ and $R(p_{\mu,\epsilon})$ are in $\Gamma$; moreover, 
$\Gamma\setminus F^{-1}_{\mu,\epsilon}(O)\subset int\,Q$ 
and $\Gamma$ is a Jordan curve with 
$ext\,\Gamma\subset B_\infty(F_{\mu,\epsilon})$. On the other hand, when 
$1<\mu\leq\mu_0'(\epsilon)$ it is used the hyperbolic nature of the fixed point 
$P_\mu$ to show that $\Gamma_b=\{(t,t):t\in 
[0,1/2]\}\cup\{(t,1-t):t\in [1/2,1]\}$ and 
$B_\infty(F_{\mu,\epsilon})=\R^2\setminus\Gamma$.

\smallskip

We complement these facts with the following:\\
{\bf Problem.}
{\em For all $\mu>1$ and $\epsilon<0$ with 
$\mu'_0(\epsilon)\leq\mu\leq\mu_2(\epsilon)$ the set of synchronized points 
$\mathcal{S}(\mu,\epsilon)$ agrees with $int\,\Gamma\cup 
F^{-1}_{\mu,\epsilon}(O)$; 
moreover, $\lim_{n\to+\infty}F^n_{\mu,\epsilon}(z)=P_\mu$ for every point
$z$ in $int\,\Gamma$, $B^c_\infty(F_{\mu,\epsilon})=cl(int\,\Gamma)$ and 
$\lim_{n\to+\infty}F^n_{\mu,\epsilon}(z)$ is either $p_{\mu,\epsilon}$ or 
$R(p_{\mu,\epsilon})$ whenever $z\in\Gamma\setminus 
F^{-1}_{\mu,\epsilon}(O)$.} 

\smallskip

As in the small strength case, we think that the region in the parameter 
space where 
$\mathcal{S}(\mu,\epsilon)=(B^c_\infty(F_{\mu,\epsilon})\setminus \Gamma)\cup 
F^{-1}_{\mu,\epsilon}(O)$ can also be increased. 

\smallskip
\noindent
2)
In relation with the existence of components of 
$B_\infty(F_{\mu,\epsilon})$ inside $int\,\Gamma$ there is a noticeable 
difference with the small strength case. When $\mu\leq\mu_1(\epsilon)$, 
so $\mu\leq 4$, it is possible that 
$B_\infty(F_{\mu,\epsilon})\cap int\,\Gamma\neq \emptyset$; indeed, this 
happens if and only if the curve $\Gamma_r$ (resp. $\Gamma_t$) passes across 
$L_1$ (resp. $L_2$). In this case every component of 
$B_\infty(F_{\mu,\epsilon})$ inside $int\,\Gamma$ is a preimage of one of the 
regions determined by either $\Gamma_r$ and $L_1$ or $\Gamma_\ell$ and $L_2$,  
as explained in item 2 of Remark \ref{rk1}; see Figure \ref{infb}. In figures 
\ref{fat3} (picture in the right side) and \ref{fat4} we show a numerical 
evidence of this fact. On the other hand, 
when $\mu>4$ it follows that $q_{\mu,\epsilon}>1$ and so the square $Q$ 
is in the topological interior of the cone $C_{\mu,\epsilon}$, 
this implies that there is no an invariante curve contained 
in $Q$ and containing $F^{-1}_{\mu,\epsilon}(O)$; indeed, 
$B^c_\infty(F_{\mu,\epsilon})$ has 
infinitely many components. 
In fact, $F^{-1}_{\mu,\epsilon}(Q)$ is disjoint union of four compact 
sets: $Q_0, Q_1, Q_2$ and $Q_3$; thus, for all $n\geq 1$ it holds that 
$F^{-n}_{\mu,\epsilon}(Q)$ is union of $4^n$ compact sets 
$Q_{i_0\cdots i_{n-1}}$ with $i_0,\cdots, i_{n-1}\in\{0,1,2,3\}$ such that:
$Q_{i_0\cdots i_{n}}\subset Q_{i_0\cdots i_{n-1}}$, 
$F_{\mu,\epsilon}(Q_{i_0\cdots i_{n}})=Q_{i_1\cdots i_{n}}$ and $z\in 
Q_{i_0\cdots i_{n}}$ if and only if $F^{j}_{\mu,\epsilon}(z)\in Q_{i_j}$ for 
every $0\leq j\leq n$. Therefore, $B^c_\infty(F_{\mu,\epsilon})$ is the union 
of the nested intersections 
$\bigcap_{n\geq 0}Q_{i_0\cdots i_n}$, where $i_n\in\{0,1,2,3\}$ for all 
$n\geq 0$. Moreover, if $\mu>4$ is large enough, then 
$B^c_\infty(F_{\mu,\epsilon})$ is a Cantor set.

\smallskip
\noindent
3) When $\mu>1$ and $\epsilon<0$ satisfy $\mu_1(\epsilon)<\mu< 4$ there is 
still the possibility of having a positively invariant curve through the points 
in $F^{-1}_{\mu,\epsilon}(O)$. 
In what follows we will describe a procedure by which we believe that such a 
curve $\Gamma$ is obtained as limit of a sequence of curves $(\Gamma_n)_{n\geq 
1}$ containing $F^{-1}_{\mu,\epsilon}(O)$. 
Take $\mu$ and $\epsilon$ as above, so $L_1$ 
(resp. $L_2$) meets $\partial Q$ in a point $q$ (resp. $R(q)$) on the segment 
$SS_1$ (resp. $S_1S_2$). The preimage under $F_{\mu,\epsilon}$ of the segments 
$OS$ and $Sq$ gives two curves: $\Gamma_1^b$ and $\Gamma_1^t$, the first one 
connects $O$ and $S$, the other one joins $S_1$ and $S_2$. In the same way, the
preimage of the segments $OS_2$ and $S_2R(q)$ is 
union of the curves: $\Gamma_1^\ell=R(\Gamma_1^b)$ and 
$\Gamma_1^r=R(\Gamma_1^t)$; recall that $R(x,y)=(y,x)$ for all $(x,y)\in\R^2$. 
Hence a closed curve containing 
$F^{-1}_{\mu,\epsilon}(O)$ is obtained: 
$\Gamma_1=\Gamma_1^b\cup\Gamma_1^t\cup\Gamma_1^\ell\cup\Gamma_1^r$. To continue 
with the construction of the sequence above
\begin{figure}[h!]
\centering
\includegraphics[scale=1]{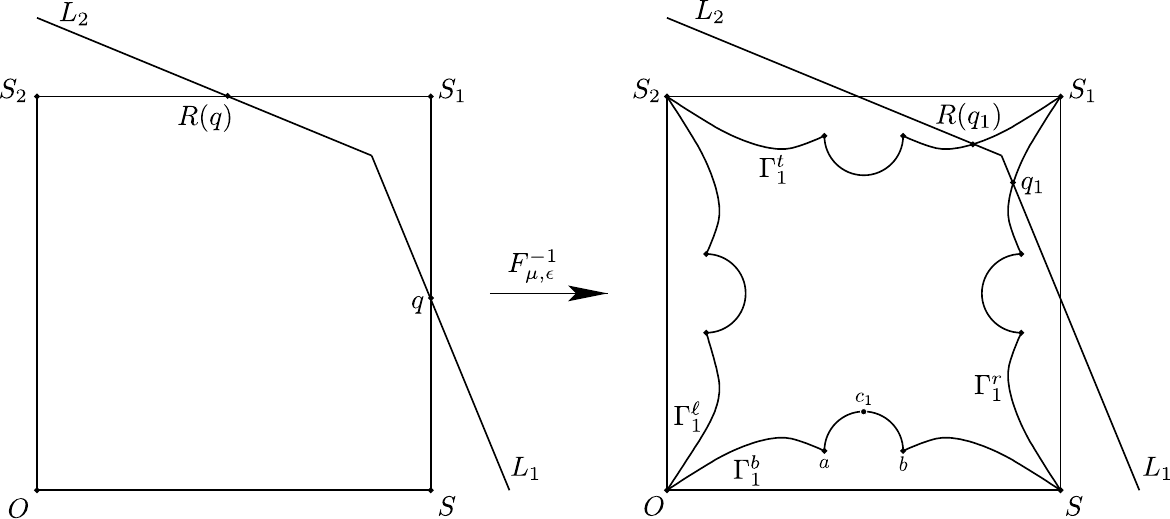}
\caption{First step in the construction of the sequence $(\Gamma_n)_{n\geq 1}$ 
which leads to a positively invariant curve through points in 
$F^{-1}_{\mu,\epsilon}(O)$.}
\label{bubbles}
\end{figure}
it is introduced the linear order $<$ on $\Gamma_1^r$: given 
$z,w\in\Gamma_1^r$, it is said that $z<w$ if and only if in the route along 
$\Gamma_1^r$, from $S$ to $S_1$, $z$ is first than $w$. In this order, let 
$q_1$ be the first point on $\Gamma_1^r$ meeting $L_1$; clearly $R(q_1)$ is the 
analagous point on $\Gamma_1^t$ from $S_2$ to $S_1$, see Figure \ref{bubbles}. 
The curve $\Gamma_2$ is obtained taking the preimage of the arcs 
$\wideparen{OSq_1}$ and $\wideparen{OS_2R(q_1)}$  on $\Gamma_1$. Next, 
repeat the same procedure, including the order above, to construct the 
remaining curves $\Gamma_n$. It should be highlighted that for each $n\geq 1$ 
the curve $\Gamma_n^b$ contains the points $a$ and $b$ of 
$F^{-1}_{\mu,\epsilon}(S)$ that are located below $\ell_2$; in addition, if 
$c_n$ is the point of 
$F^{-1}_{\mu,\epsilon}(q_{n-1})$ ($q_0=q$) in $\Gamma_n^b$, then  
the arc $\wideparen{ac_nb}$ is the part of the preimage of the piece on 
$\Gamma_{n-1}^r$ from $S$ to $q_{n-1}$ ($\Gamma_0^r$ is the segment $SS_1$ on 
$\partial Q$); that arc is a kind of bubble that connects the two arcs in 
$\Gamma_n^b$ produced by means of the preimage of $\Gamma_{n-1}^b$; see Figure 
\ref{bubbles} for the case $n=1$.

\smallskip

A couple of important remarks about the sequence $(\Gamma_n)_{n\geq 1}$ and its
limit curve $\Gamma$ are necessary. First, for $n$ large 
enough the curve $\Gamma_n^b$, even $\Gamma^b$, may not be the graph of a 
function; this loss is essentially due to the fact that the origin is a 
hyperbolic repeller and the expansion on the antidiagonal is greater than the 
expansion on the diagonal, which implies that the slope of 
$\Gamma_n^b$ at $S$ is asymptotic to $-1$. Thus, $\Gamma_n^b$ and $\Gamma_n^r$ 
look tangent to the antidiagonal at $S$ when $n$ is large enough.
\begin{figure}[h!]
\centering
\includegraphics[scale=1]{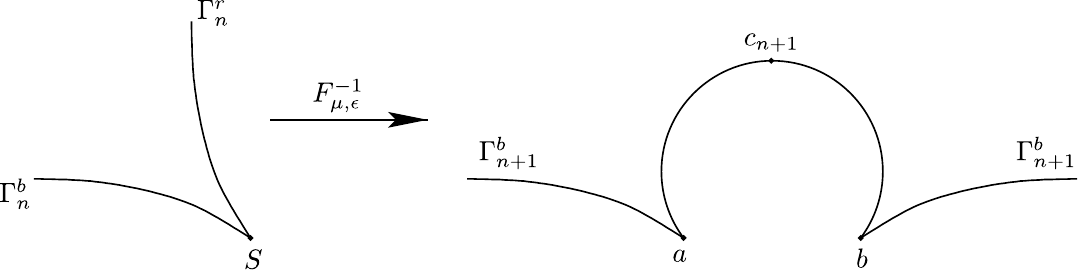}
\caption{Asymptotic tangency of curves $\Gamma_n^b$ and $\Gamma_n^r$ at $S$ 
implies the breaking of the graph condition for the curve $\Gamma_{n+1}^b$ 
near the points $a$ and $b$.}
\label{nofunction}
\end{figure}
This configuration is propagated by preimages, so in a neighborhood of the 
point $a$, and also at $b$, the vertical line test indicates that 
$\Gamma_n^b$ 
is not the graph of a function when $n$ is large enough. Figure 
\ref{nofunction} illustrates how the curve $\Gamma_n^b$ loses the graph 
nature when $n$ is large enough, see also Figure \ref{nodif}. The second remark 
is related to the possibility that 
$ext\,\Gamma\cap B^c_\infty(F_{\mu,\epsilon})\neq\emptyset$. To be able to 
argue this possible phenomenon, suppose that 
for some $n\geq 1$ there are points $r_1,r_2\in\Gamma_n^r$ such that, in the 
order considered in $\Gamma_n^r$, it holds that  $q_n<r_1<r_2$ and the arc in 
$\Gamma_n^r$ containing $r_1$ and $r_2$ is contained in the 
cone $C_{\mu,\epsilon}$; so, the preimage of this arc (which is not considered 
to obtain $\Gamma_{n+1}$) is union of two closed curves contained in 
$ext\,\Gamma_{n+1}$; see Figure \ref{outside} where this possibility is grossly 
illustrated.
\begin{figure}[h!]
\centering
\includegraphics[scale=1]{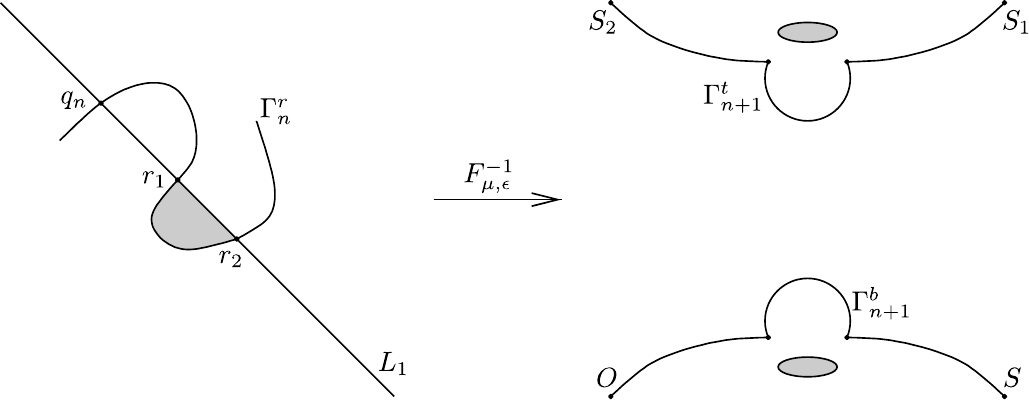}
\caption{An outline of the configuration on $\Gamma_n$ that provides 
bounded orbits in $ext\,\Gamma_{n+1}$.}
\label{outside}
\end{figure}
Observe that if this arrangement is also on the limit curve $\Gamma$, then 
the preimage of an arc as described for $\Gamma_n^r$ but now on $\Gamma^r$ 
produces components of $B_\infty(F_{\mu,\epsilon})$ in $ext\,\Gamma$ whose 
boundaries are points with omega limit set contained in $\Gamma$.
\end{rk}

\section{Numerical evidences}\label{compu}
In this section we will report part of the results of computational simulations 
that we have done to try to understand a little more the complicated dynamics of
$F_{\mu,\epsilon}$. Indeed, we display some figures 
generated by forward or backward iterations of certain points in order to try 
to obtain information about $B_\infty(F_{\mu,\epsilon})$, its  boundary and 
even other dynamic behaviors that are not theoretically described. 
Figures 
shown below were generated using Fractint, a freeware fractal generator  
(https://www.fractint.org/). By means of this platform we illustrate 
several results discussed in \cite{rsv} 
and \cite{rrv}, and some other possible dynamic phenomena that by the 
complicated 
nature and the abscense of appropriate theoretical tools, they have not 
been addressed with adequate mathematical rigor. Indeed, we hope that based on 
these computational simulations one can follow the heuristic principle for 
further analytical examination that allows to explain phenomena 
such as Hopf bifurcations or existence of fat
attractors, which, due to the numerical evidences, could be displayed in the 
dynamic of the coupled logistic map \eqref{quad}.

\smallskip

In this computational platform we introduce a program to plot 
forward and backward iterations of points under $F_{\mu,\epsilon}$, it is also 
plotted preimages of the circle $C$ and the boundary of the square $Q$. 
Some comments are necessary before showing a collection of figures 
produced with Fractint. When the backward iteration of a point is done, the 
points plotted outline a figure in which two zones appear: a dark one 
(preimages of the point) containing part of the complement of 
$B_\infty(F_{\mu,\epsilon})$ and a white zone that does not show preimages of 
the point by one of the following reasons: either the number of iterated is 
insufficient or such white zones are part of the basin of attractors, one of 
them the attractor at $\infty$ and possibly other attractors contained in the 
bounded component of the boundary of the immediate basin of $\infty$.

\smallskip

Our program is concentrated in five subroutines:

\smallskip
\noindent
{\bf Option $\boldsymbol 0$.} Given an integer $n\geq 0$ and $(x,y)\in\R^2$, 
every point $F^k_{\mu,\epsilon}(x,y)$ with $0\leq k\leq n$ is plotted. With 
this option possible attractors should be detected.

\smallskip
\noindent
{\bf Option $\boldsymbol 1$.} 
For $n\geq 0$ and $(x,y)\in\R^2$ all the possible points in 
$F^{-n}_{\mu,\epsilon}(x,y)$ are plotted; recall 
that $(x,y)$ has preimage if and only if $(x,y)\in C_{\mu,\epsilon}$. This 
tool allows to get an idea of the invariant curve $\Gamma$ in theorems 
\ref{a}, \ref{b} and remarks \ref{rk1}, \ref{rk2}.

\smallskip
\noindent
{\bf Option $\boldsymbol 2$.}  This option is a combination of the previous 
ones. Given a point $(x,y)$ and nonnegative integers $n_1,n_2$ with 
$n_1+n_2>0$, for each $0\leq j\leq n_1$ Option 1 is applied to 
$F^j_{\mu,\epsilon}(z)$ with $n=n_2$.

\smallskip
\noindent
{\bf Option $\boldsymbol 3$.} It is a subroutine to plot preimages of the 
circle $C$; it is a tool to visualize approximations of the 
curve $\Gamma$ in the small strength case.

\smallskip
\noindent
{\bf Option $\boldsymbol 4$.} This is similar to Option 3 but $C$ 
is replaced by the boundary of square $Q$. It is useful to have 
information about $\Gamma$ in the large strength case.

\smallskip

In each 
of the figures shown in paragraphs 
\ref{small} and \ref{large} are notified the values of $\mu$ and $\epsilon$, 
in some cases option used is also informed.

\subsection{Small strength case}\label{small}
Here we discuss some scenarios of how the curve $\Gamma$ and its preimages are 
arranged. With the help of computational simulations we also show a colections 
of figures evidencing several kind of attracting sets, some of which are 
related with possible Hopf bifurcations.

\subsubsection{\bf On the synchronization case}
\begin{figure}[h!]
\centering
\includegraphics[scale=1.2]{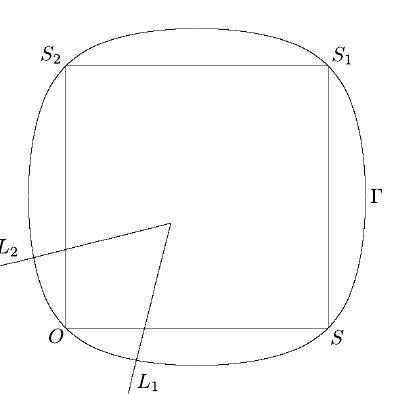}
\hspace{1cm}
\includegraphics[scale=1.2]{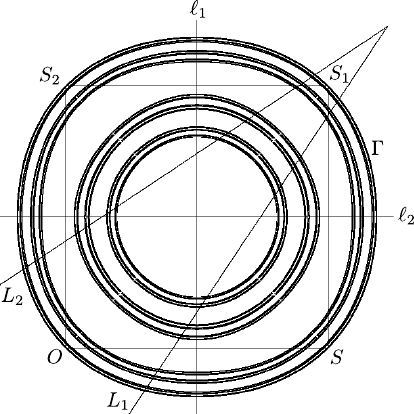}
\caption{In these pictures are shown large backward iterations of the 
circle $C$. In the 
left one $\mu=1.6$ and $\epsilon=0.2$, for the picture on the right side 
$\mu=4.9$ and $\epsilon=0.4$. The indicated curve $\Gamma$ is actually 
an approximation to the Jordan curve in Theorem \ref{a}; that  
approximation in the right side is just the external border of that picture.}
\label{figura12}
\end{figure}

Consider the region in the parameter space where it has been proved 
that synchronization occurs (see item \ref{a-c} of Theorem \ref{a}), that is  
$1<\mu\leq\max\{\mu'(\epsilon),\mu_0(\epsilon)\}$. For all $\mu$ and $\epsilon$ 
in this sector it is clear that $B_\infty(F_{\mu,\epsilon})=ext\,\Gamma$ and 
$\partial B_\infty(F_{\mu,\epsilon})=\Gamma$; in addition, every point inside 
$int\,\Gamma$ has its omega limit set contained in the segment $OS_1$. However 
the fashion in which preimages of $\Gamma$ are arranged is not uniform in that 
parameter region. In fact, if $\mu<4$, then $\Gamma$ is completely invariant: 
$F^{-1}_{\mu,\epsilon}(\Gamma)=\Gamma$; the picture in the left side in Figure 
\ref{figura12} shows the manner in which $\Gamma$ and its preimage appear in 
this case. When $\mu=4$, the countable set $\bigcup_{n\geq 
0}F^{-n}_{\mu,\epsilon}(\Gamma)$ has a non-trivial part inside $int\,\Gamma$; in 
fact, the set of points $q$ in the segment $OS_1$ such that 
$F^{n}_{\mu,\epsilon}(q)=S_1$ for some $n\geq 0$ is dense in that segment. Now, 
when $\mu>4$ every component of the preimages of $\Gamma$ is a continuum 
(compact and connected) set. In what follows we discuss about this fact; before 
a possible configuration of 
$\bigcup_{n\geq 0}F^{-n}_{\mu,\epsilon}(\Gamma)$ is shown in right picture 
of Figure \ref{figura12}. 
Recall that if $\mu>4$, then 
$B_\infty(F_{\mu,\epsilon})$ has infinitely many components in 
$int\,\Gamma$; see first part of item 2 in Remark \ref{rk1}. These 
components are arranged generically in two special ways; to explain 
them consider the set $P(O)$ of all the
preimages of $O$ on the diagonal, also consider the Cantor set $K_\mu$ defined 
dynamically by the logistic map $f_\mu$; clearly $P(O)$ is a proper subset of 
$K_\mu$. Let 
$T$ be the triangular region introduced at the beginning of item 2 in Remark 
\ref{rk1}. So, $F^{-1}_{\mu,\epsilon}(T)$ is a topological disk whose interior 
is contained in $B_\infty(F_{\mu,\epsilon})$ and its
boundary 
is the Jordan curve obtained as preimage of the arc $\wideparen{p_1p_2}$ on 
$\Gamma$ as defined in Remark \ref{rk1}, see Figure \ref{infa}. Recall 
that for a 
simple curve inside $C_{\mu,\epsilon}$ with extreme points at $L_1$ and $L_2$, 
its preimage is always a topological circle (a circle, for short) surrounding 
$c=(1/2,1/2)$; this preimage is
generically either in the interior of $C_{\mu,\epsilon}$ or its intersection 
with that cone is union of two disjoint curves with extreme points at $L_1$ and 
$L_2$. From now on we say that a curve {\em goes 
through} $L_1$ and $L_2$ if it is contained inside $C_{\mu,\epsilon}$ with
extreme points in $L_1$ and $L_2$.

\smallskip
\noindent
{\bf First way:} {\em The boundary of $F^{-1}_{\mu,\epsilon}(T)$ 
goes through $L_1$ and $L_2$.}\\
In this instance, 
the preimage of $F^{-1}_{\mu,\epsilon}(T)$ is a closed topological annulus 
(an annulus, for short) such that: its interior 
is part of $B_\infty(F_{\mu,\epsilon})$ and each circle of its boundary
cuts $P(O)\setminus \{O,S_1\}$ in two points and goes 
through $L_1$ and $L_2$. 
Repeating this argument recursively one can describe the components of 
$B_\infty(F_{\mu,\epsilon})$ inside $int\,\Gamma$. The following facts 
are clears:
\begin{enumerate}[(a)]
\item
For every $n\geq 2$, $F^{-n}_{\mu,\epsilon}(T)$ is disjoint union of $2^{n-2}$ 
annuli whose interiors is in $B_\infty(F_{\mu,\epsilon})$; as above, each 
circle in the boundary of every 
annulus in $F^{-n}_{\mu,\epsilon}(T)$ goes through $L_1$ and $L_2$, and also 
cuts $P(O)\setminus \{O,S_1\}$ exactly in 
two points.
\item
$P(O)=\{O,S_1\}\cup\bigcup_{n\geq 1}(\partial 
F^{-n}_{\mu,\epsilon}(T)\cap\Delta)$, where $\Delta$ is, as above, the 
diagonal.
\end{enumerate}
In addition, since each point $q\in K_\mu\setminus P(O)$ also has a stable 
manifold $W^s(q)$, with arguments like $\lambda$-lemma (see 
\cite{palis}) one can prove that the component of $W^s(q)$ containing $q$ is a 
circle accumulated by the boundaries of the annuli above described. 
Consequently the complement of 
$B_\infty(F_{\mu,\epsilon})$ is a Cantor set of circles; in fact, every pair of 
points in $K_\mu$ with same image under $F_{\mu,\epsilon}$ are joined by such a 
circle. Observe that this Cantor set of circles is just the closure of 
$\bigcup_{n\geq 0}F_{\mu,\epsilon}^{-1}(\Gamma)$. Right picture in 
Figure \ref{figura12} is an illustration of that phenomenon.

\smallskip
\noindent
{\bf Second way:} {\em $F^{-1}_{\mu,\epsilon}(T)$ is contained into 
the interior of $C_{\mu,\epsilon}$.}\\
This case is really of higher complexity than the previous one. To facilitate
our unfinished attempt to describe how the components of 
$B_\infty(F_{\mu,\epsilon})\cap int\,\Gamma$ are arranged we introduce the 
following concept. An annulus is said to be {\em large} if each 
circle in its boundary goes through $L_1$ and $L_2$.

\smallskip

\begin{figure}[h!]
\centering
\includegraphics[scale=1]{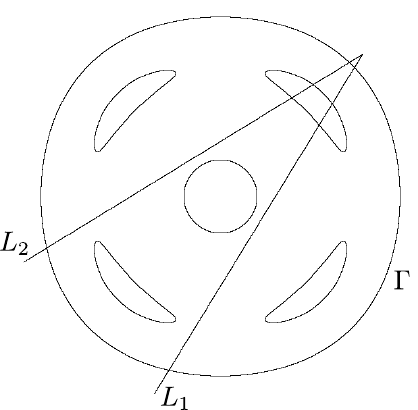}
\hspace{1cm}
\includegraphics[scale=1]{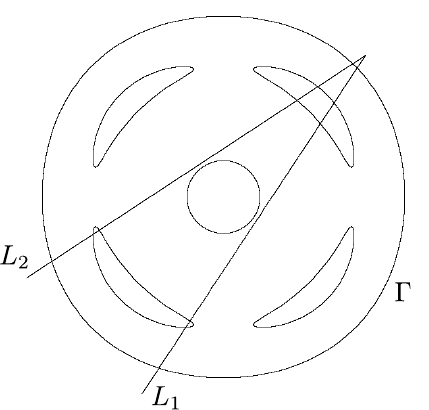}

\includegraphics[scale=1.33]{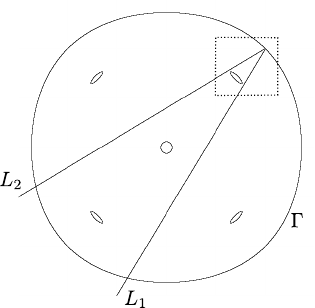}
\hspace{1cm}
\includegraphics[scale=1.33]{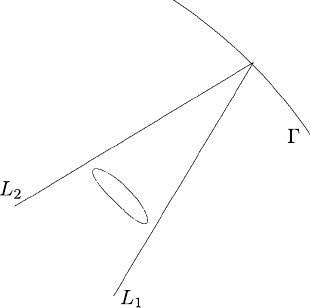}
\caption{The pictures of the top and the picture in left side of the bottom
illustrate the generic configurations of the 
preimage of $F^{-1}_{\mu,\epsilon}(T)$ (central disk). In left top picture
$\epsilon=0.38$ and $\mu=4.16$, in right one $\epsilon=0.397$ and $\mu=4.16$;
in left bottom picture
$\epsilon=0.3755$ and $\mu=4.007$. The picture in the right side of the bottom
picture is a magnification of the dotted square. Observe that in all of these 
computational simulations
$\mu<\mu_0(\epsilon)$.}
\label{figura14}
\end{figure}

Due to the location of the disk $F^{-1}_{\mu,\epsilon}(T)$, two 
facts follow:

\smallskip
\noindent
$\bullet$
The preimage of $F^{-1}_{\mu,\epsilon}(T)$ is disjoin union of four disks, 
the interior of each disk is part of 
$B_\infty(F_{\mu,\epsilon})$; further, 
two of them intersect $P(O)\setminus \{O,S_1\}$ and the other two have no 
preimage. Generically these four disk are displayed as illustrated in Figure 
\ref{figura14}. All pictures in this figure were generated using Option 3.

\begin{figure}[h!]
\centering
\includegraphics[scale=1.2]{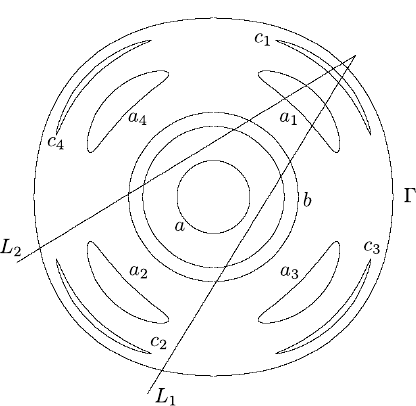}
\hspace{1cm}
\includegraphics[scale=1.2]{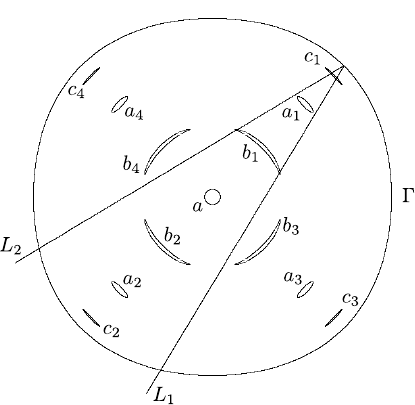}
\caption{In both picture, the disks denoted by $a_1,a_2,a_3$ and $a_4$ are the 
preimage of $F^{-1}_{\mu,\epsilon}(T)$ (central disk). The annulus $b$ in left 
picture is the preimage of the disk $a_1$, and the disks $b_1,b_2,b_3$ and 
$b_4$ constitute the preimage of $a_1$ in picture on the right side.}
\label{figura15}
\end{figure}

\begin{figure}[h!]
\centering
\includegraphics[scale=1]{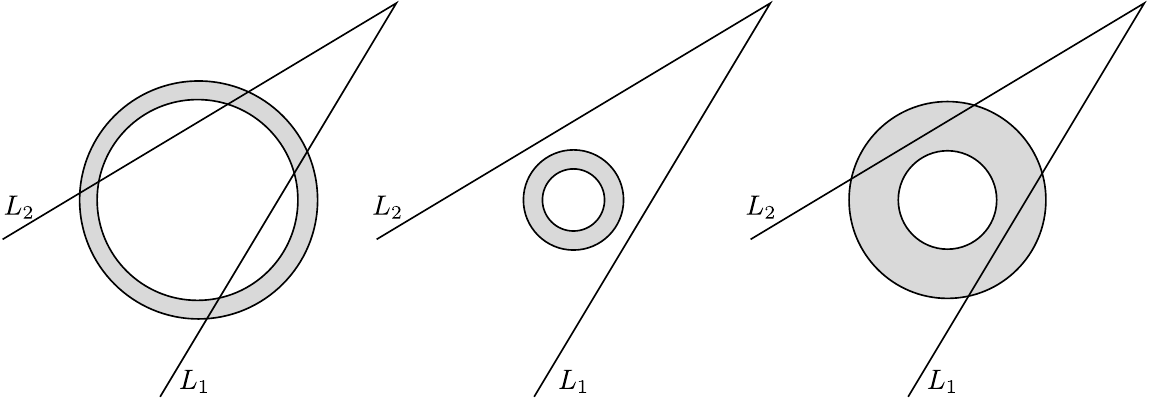}
\caption{The three generic possibilities for the preimage of a disk whose 
boundary in $C_{\mu,\epsilon}$ goes through $L_1$ and $L_2$.}
\label{figura16}
\end{figure}

\noindent
$\bullet$
Since $O$ is a hyperbolic saddle, from $\lambda$-lemma it follows that if $A$ 
is a disck in $F^{-2}_{\mu,\epsilon}(T)$ intersecting $\Delta$, then there 
exist $k\geq 1$ and a component $B$ of $F^{-k}_{\mu,\epsilon}(A)$ such that $B$ 
is a disk whose 
boundary goes through $L_1$ and $L_2$; pictures in Figure \ref{figura15} show 
this situation for 
$\epsilon=0.38$ and $\mu=4.16$ (picture in the left side), and 
$\epsilon=0.3755$ 
and $\mu=4.007$ (picture in the right side).
The preimage $F^{-1}_{\mu,\epsilon}(B)$ is not necessarily a large annulus; 
indeed, it is generically displayed in one of the three ways as shown in Figure 
\ref{figura16}; right there the left annulus is the only large; central and 
right annuli will be called {\em small} and {\em singular} annulus, 
respectively. The preimages of these three type of annuli are displayed of 
three different manners, from
left to right: two annuli, at least one of them large; four annuli, none large; 
and one annulus minus four disks. An interesting problem is the following:

\smallskip
\noindent
{\bf Problem.} 
{\em In the second way setting, give a detailed description of all 
manners in which the basin of attraction of $\infty$ in 
$int\,\Gamma$ is displayed.} 

\smallskip

\begin{figure}[h!]
\centering
\includegraphics[scale=1]{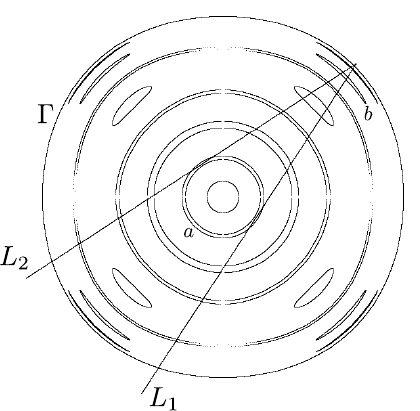}
\hspace{1cm}
\includegraphics[scale=1]{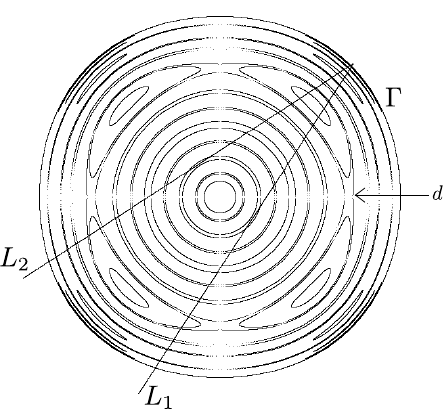}

\includegraphics[scale=1]{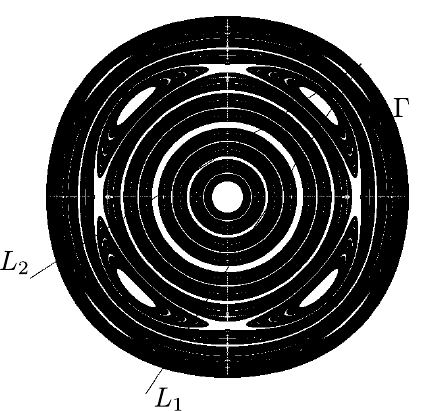}
\caption{Computational illustration through which the 
existence of singular annuli is evidenced. White sectors at figure bottom 
illustrates some components of $B_\infty(F_{\mu,\epsilon})$.}
\label{figura17}
\end{figure}

We advance that among the components of $B_\infty(F_{\mu,\epsilon})\cap 
int\,\Gamma$ always there are large and small annuli.
Given a large annulus in $B_\infty(F_{\mu,\epsilon})\cap int\,\Gamma$, the 
preimage of 
its components in $C_{\mu,\epsilon}$ are two annuli; one of them is large 
and its components in $C_{\mu,\epsilon}$ are closer to $O$ and $S_1$ than the 
previous ones. The preimage of the other component may be 
large, small or singular, it is closer to 
$F^{-1}_{\mu,\epsilon}(T)$ than the given annulus. So, from $\lambda$-lemma and 
continuity it follows that $\Gamma$ is accumulated by large annuli; in fact, it 
is accumulated by preimages of closer component to $O$ of large annuli. 
Consequently the stable manifold of $O$ is accumulated by preimages of large 
annuli. It is also clear that the boundary of $F^{-1}_{\mu,\epsilon}(T)$ is 
accumulated by small annuli.
On the other hand, observe that the number of disks in 
$\bigcup_{n\geq 1}F^{-n}_{\mu,\epsilon}(T)$ is finite. This is consequence of 
the following facts: those disks in $\bigcup_{n\geq 1}F^{-n}_{\mu,\epsilon}(T)$ 
that do not cut $\Delta$ have no preimage, and among those who cut the 
diagonal, only two of them are such that the boundary of each one goes through 
$L_1$ and $L_2$; one of these two disks is near $O$ and the other one is near 
$S_1$. The preimage of the first of these disks is a large annulus (as it is 
known); the preimage of the second one may be a large, small or singular 
annulus. Suppose it is singular, hence any other annulus surrounding 
$F^{-1}_{\mu,\epsilon}(T)$ must be either large or a small annulus. 
The same thing happens if a singular annulus surrounds any other 
of the disks above; thus, there can only be  a finite number of singular 
annuli. Figure \ref{figura17} shows a numerical evidence of the existence of a 
singular annulus when $\mu=4.03$ and $\epsilon=0.394$. Pictures at the top of 
this figure were obtained with Option 3; in the left picture observe that the 
annulus denoted by $a$ is a singular annulus and it is the preimage of the disk 
indicated with $b$. 
Sector $d$ (at the right) is the preimage of $a$. Picture at the bottom was 
generated by means of Option 2 with starting point $(0.1,0)$, region in black 
corresponds to huge numbers of forward and backward iterations of this point.

\subsubsection{\bf Attractors inside $\boldsymbol{int\,\Gamma}$}
In a wide variety of dynamical systems there are several coexisting attractors, 
many of them easy to be detected both numerically and theoretically. 
To fully understand the dynamics of these systems one should identify all the 
possible attractors and their basins of attractions. At this point the 
computational simulation is actually an important tool, it allows the location 
of possible attractors and their basins. 

\smallskip

To specify terms we recall that if 
$T$ is a continuous self-map on a metric space $X$, a compact subset $A$ of $X$ 
is called {\em attracting set} for $T$ if $T(A)\subset A$ ($A$ is invariant) 
and there exists a neighborhood $V$ of $A$ such that for any neighborhood $U$ 
of $A$ there exists an integer $n_U$ such that for any $n>n_U$, $T^n(V)\subset 
U$. Observe that for any point $x\in V$ it holds that 
$T^n(x)\to A$ when $n\to +\infty$. The {\em basin of attraction} of an 
attracting set $A$ is the set of points whose forward orbit accumulates on $A$; 
clearly the basin of attraction of $A$ is the equal to $\bigcup_{n\geq 
0}T^{-n}(V)$. An {\em attractor} for $T$ is any attracting 
set containing a dense forward orbit. We remark that an 
attracting set may contain one or several attractors. We say that an attracting 
set is a {\em fat attractor} if it has positive Lebesque measure. The term fat 
attractor is also used in other settings, see for example \cite{lopes}, 
\cite{mora} and \cite{tsujii}. We refer to \cite{milnor,milnorb} for further 
discussion on the definition of attractors.

\smallskip

In this paragraph we show some numerical evidences of the existence of fat 
attractors for $F_{\mu,\epsilon}$. First, observe that for  
$1<\mu\leq\max\{\mu'(\epsilon),\mu_0(\epsilon)\}$, 
there is nothing relevant to say about the forward 
asymptotic behavior of the orbits starting in $cl(int\,\Gamma)$:  
the omega limit set of all $z\in cl(int\,\Gamma)$ is either contained in 
$\{R(p_{\mu,\epsilon}),p_{\mu,\epsilon}\}$ or in the segment $OS_1$. 
Therefore we will consider the following two cases for the parameter regions: 
$0<\epsilon<3/8$ and $\mu'(\epsilon)<\mu\leq 4$, and 
$\max\{4,\mu_0(\epsilon)\}<\mu\leq \mu_1(\epsilon)$; see Figure \ref{3curvas1}.

\smallskip
\noindent
{\bf a) Case 1.} $0<\epsilon<3/8$ and $\mu'(\epsilon)<\mu\leq 4$.\\ 
From Theorem \ref{a} it is clear that $cl(int\,\Gamma)= 
B_\infty^c(F_{\mu,\epsilon})$ and $ext\,\Gamma=B_\infty(F_{\mu,\epsilon})$.  
We believe, according to computational simulations, that in the parameter 
region above there are sets with positive measure for which 
$F_{\mu,\epsilon}$ exhibits fat attractors. A theoretical 
proof of this statement seems very hard, so we only show 
several images as numerical evidences of this occurrence.

\smallskip

\begin{figure}
\centering
\includegraphics[scale=1]{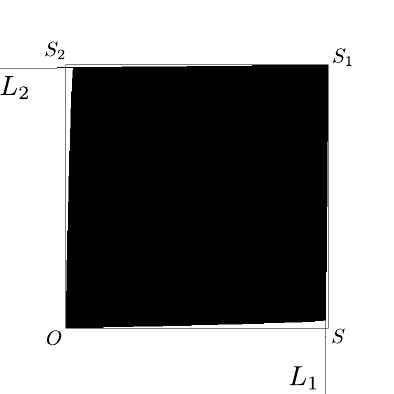}
\includegraphics[scale=1]{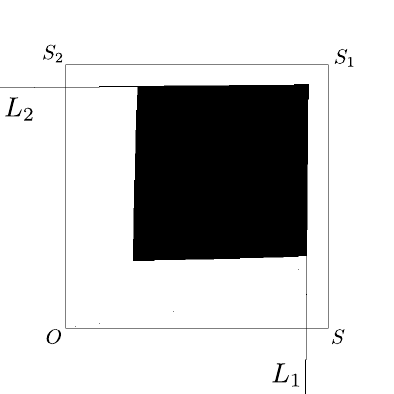}
\includegraphics[scale=1]{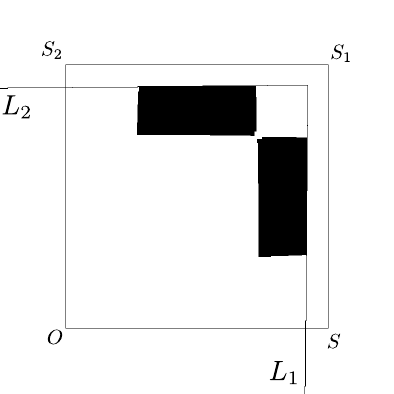}

\includegraphics[scale=1]{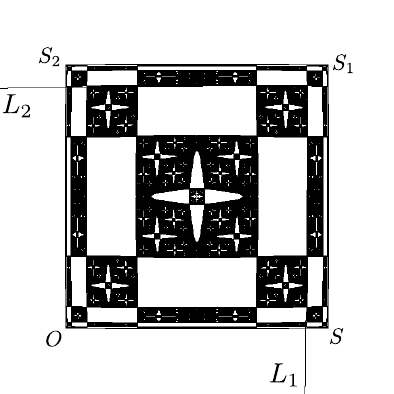}
\hspace{1cm}
\includegraphics[scale=1]{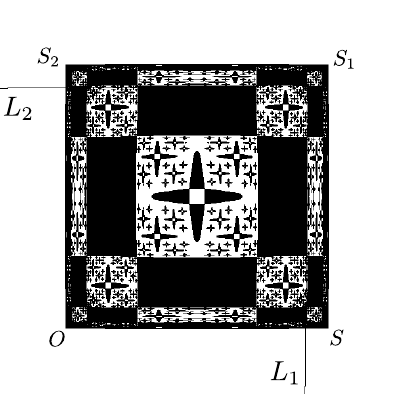}
\caption{The parameters for the picture in the left side are $\mu=4$ and 
$\epsilon=0.01$, for the central picture $\mu=3.7$ and $\epsilon=0.01$ and for 
the right picture $\mu=3.694$ and $\epsilon=0.01$; this last fat attractor is 
$2$-periodic. Pictures in the bottom were obtained through Option 2 to 
visualize an approximation of the basins of attraction of the attracting set in 
the segment $OS_1$ (left picture) and the fat attractor when $\mu=3.694$ and 
$\epsilon=0.01$ (right picture).}
\label{fat1}
\end{figure}

The first attempt to detect fat attractors is for parameter values close to 
$\epsilon=0$ and $\mu=4$, this is because for these values of $\epsilon$ and 
$\mu$ almost every point inside $int\,\Gamma$ has dense orbit in the square 
$Q$. So, for $\mu\leq 4$ and $\epsilon>0$ close enough to $4$ and $0$ 
respectively, the map $F_{\mu,\epsilon}$ should also exhibit fat attractors.
Figure \ref{fat1} shows fat attractors (black regions) when 
$\epsilon=0.01$ and $\mu$ takes 
three different values: $4, 3.7$ and $3.694$
These regions were obtained by plotting a huge number $(2\times 10^8)$ of 
forward iterated of points in the square $Q$. For that parameters, the 
black regions in that figure illustrate the only fat attractors of
$F_{\mu,\epsilon}$. However there may be other attractors, for example when 
$\epsilon=0.01$ and $\mu=3.694$ the map 
$F_{\mu,\epsilon}$ has an attractor on the diagonal. We justify this statement 
through computational simulation. Recall that if one takes any point on the 
segment $OS_1$, then its forward orbit remains in this segment and accumulates 
on the attractor of the logistic map $f_\mu$; it is obvious that the forward 
orbit of every point in any preimage set of that point also accumulates on that 
attractor. Therefore Option 2 is a nice tool to detect the basin of possible 
attractors. We have done this computational simulation with $\epsilon=0.01$ and 
$\mu=3.694$, the results are the pictures in the bottom of Figure \ref{fat1}.

\smallskip

\begin{figure}[h!]
\centering
\includegraphics[scale=1]{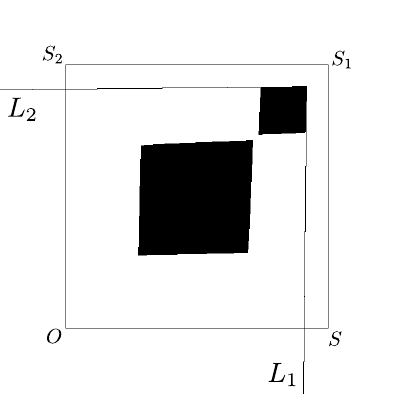}
\hspace{1cm}
\includegraphics[scale=1]{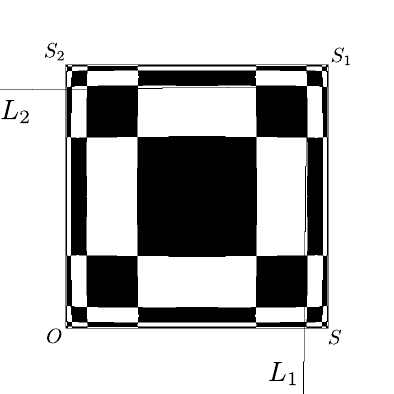}

\includegraphics[scale=1]{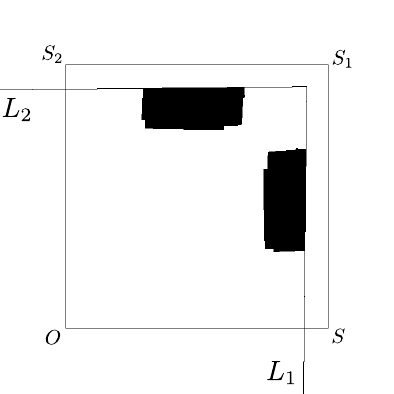}
\hspace{1cm}
\includegraphics[scale=1]{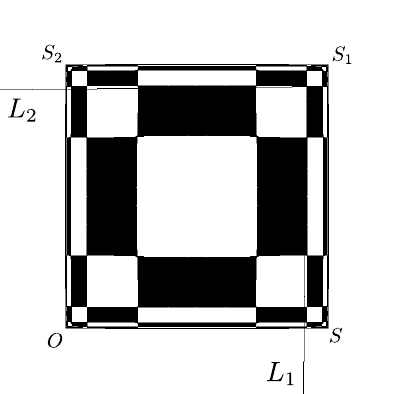}
\caption{Two different fat attractors for the same map $F_{\mu,\epsilon}$, here 
$\mu=3.67$ and $\epsilon=0.01$.}
\label{fat2}
\end{figure}

For other values ​​of the parameters $\mu$ and $\epsilon$ there may be more 
than one fat attractor, this is numerically evidenced in the pictures 
of Figure \ref{fat2}.
The two fat attractors in that figure are located in the 
left side (bottom and top), pictures in the right side are approximations 
of the basin of attraction (black sectors) of the corresponding attractor to 
its left side. Such basins of attraction were obtained by using Option 2 as 
explained above. Observe that these attracting sets are periodic with 
period 2 and their basins look complementary, so there are no more attractors 
when $\mu=3.67$ and $\epsilon=0.01$.

\smallskip

\begin{figure}
\centering
\begin{tikzpicture}
\pgftext{\includegraphics[scale=1]{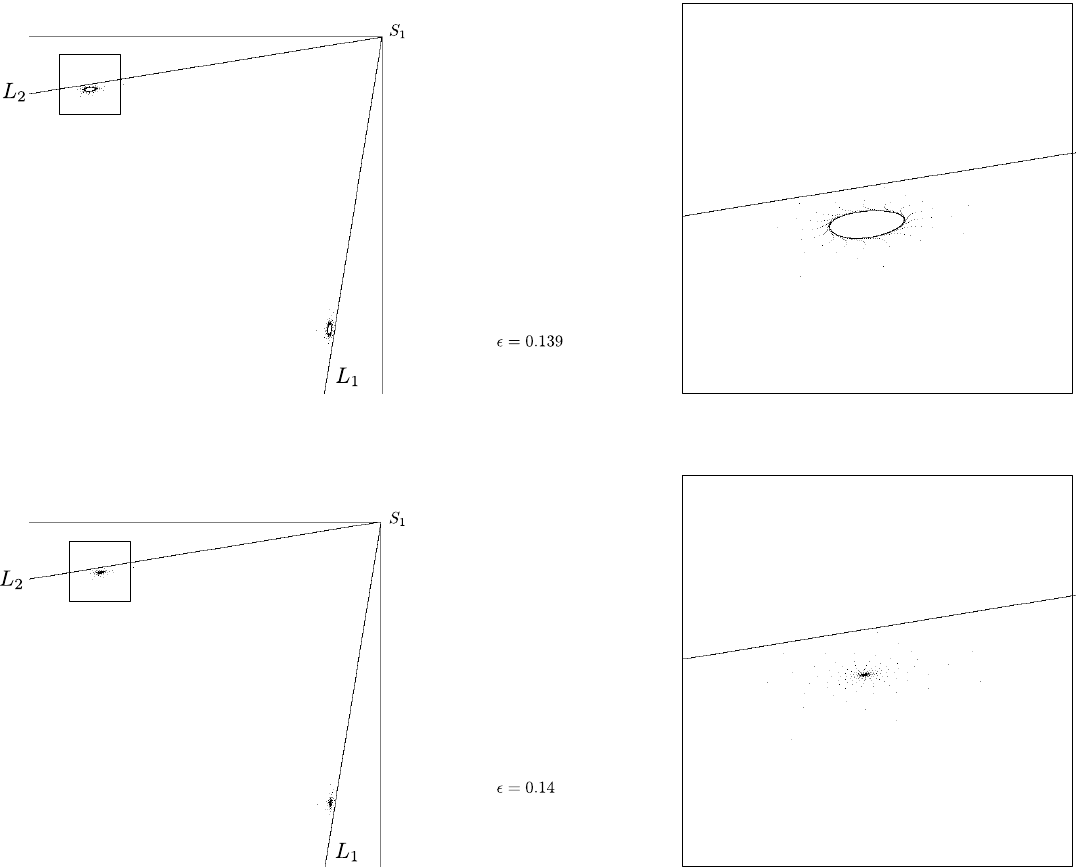}}
\draw[black!50, very thin,->] (-4.2,3.4) -- (-0.5,3.4) -- (-1,2.9) -- 
(1.4,2.9);
\draw[black!50, very thin,->] (-4.1,-1.5) -- (-0.5,-1.5) -- (-1,-2) -- 
(1.4,-2);
\end{tikzpicture}
\caption{Numerical evidences of the evolution of the fat attractors going 
through a Hopf bifurcation when $\epsilon$ increases from $0$ and $\mu=4$.}
\label{mu4small4}
\end{figure}

To finish with this case we will report some numerical evidences revealing
\begin{figure}[h!]
\centering
\includegraphics[scale=1]{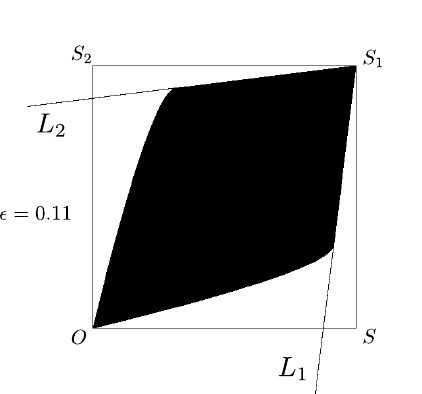}\hspace{1cm}
\includegraphics[scale=1]{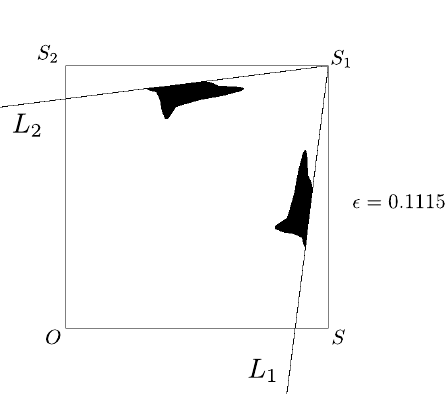}

\vspace{.5cm}
\begin{tikzpicture}
\pgftext{\includegraphics[scale=1]{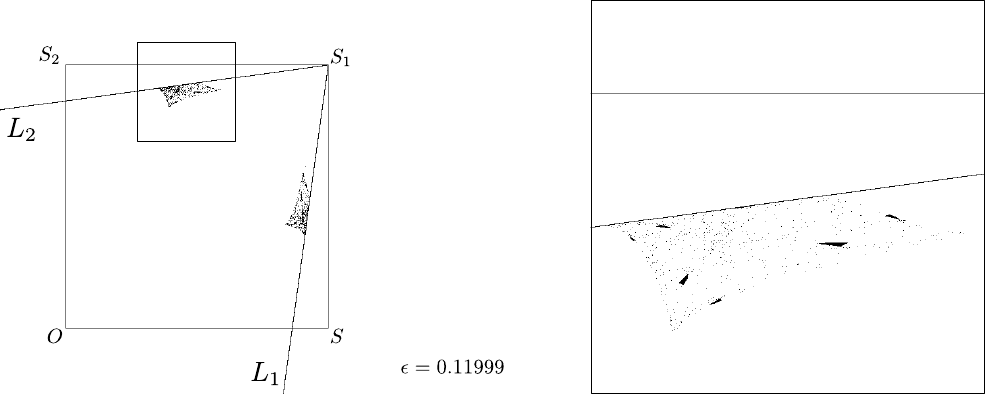}}
\draw[black!50, very thin,->] (-2.55,1) -- (-0.5,1) -- (-1,0.5) -- (0.9,0.5);
\end{tikzpicture}

\includegraphics[scale=1]{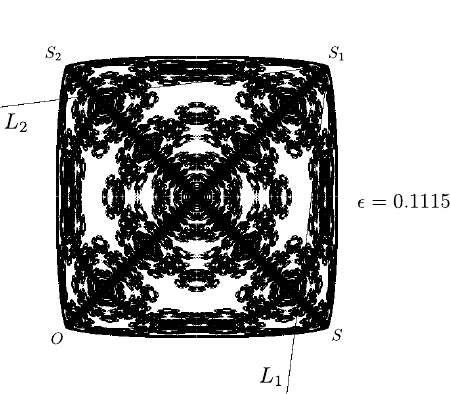}
\caption{This figure marks the beginning of how fat attractors change when 
$\mu=4$ and $\epsilon$ grows from $0$. Pictures at the top and in the middle 
were obtained with Option $0$, while the image at the bottom was generated 
with Option 1 applied to the origin, it also can be obtained with Option 2 
applied to any point in the segment $OS_1$; sectors in white constitute 
what seems to 
be the basin of attraction of the $2$-periodic fat attractor at the top.}
\label{mu4small}
\end{figure}
not only the existence of fat attractors when $\mu=4$ but also the 
dynamics richness in their evolution. We recall that the logistic map $f_4$ has 
the interval $[0,1]$ as a completely invariant set and its dynamics is chaotic 
in Devaney's sense; see \cite{devaney}. This fact is perhaps not exactly 
relevant for the following discussion, but it does provides information about 
the dynamics of $F_{4,\epsilon}$ restricted to the diagonal.
In Figure \ref{fat1} we 
have shown a such attracting set for $\epsilon=0.01$, this kind of attracting 
sets persist when $\epsilon$ increases at least until $\epsilon=0.11$; see the 
fat attractor located in the left 
side of the top of the Figure \ref{mu4small}. The 
$2$-periodic fat attractor to its right 
side appears when $0.11<\epsilon\leq 0.1115$, in this case $\epsilon$ is 
exactly $0.1115$; it is a break of the continuum of fat attractors when 
$\epsilon>0$ grows from $0$. Clearly this attracting set does 
not contain the segment $OS_1$, which we believe is an 
attractor. 
The reason of our belief is based on the fact that by plotting 
a huge number of preimages of the origin, indeed using Option 2 to any point 
in $OS_1$, the result is a set of points whose complement should be the basin 
of attraction of that $2$-periodic fat attractor; picture at the bottom of 
Figure \ref{mu4small} shows that plotting set when $\epsilon=0.1115$.
To produce the picture in the middle of Figure \ref{mu4small}, that is when 
$\epsilon=0.11999$, we have used Option 0; the result is an 
apparently diffuse set of points. However, when one does a magnification in on 
one of its two parts, small fat regions can be observed; just with these 
regions a 
periodic fat attractor is constituted: the $2$-periodic fat attractor transits 
by an explosion in its parts to produce a higher order periodic fat attractor, 
see the zoom of the box in the image mentioned above. After that explosion the 
pieces of the new periodic fat attractor seem to come together to form a new 
type of attracting set, Figure \ref{mu4small3} shows the state of that evolution 
when $\epsilon=0.125$.
\begin{figure}
\centering
\begin{tikzpicture}
\pgftext{\includegraphics[scale=1]{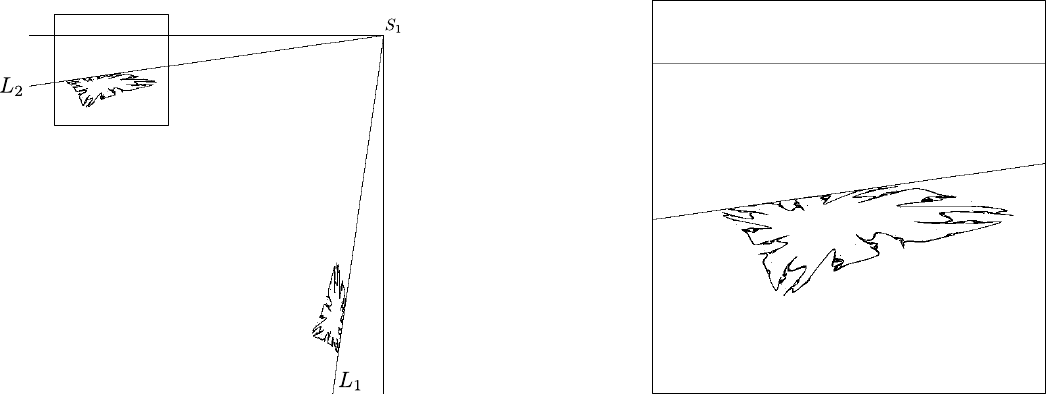}}
\draw[black!50, very thin,->] (-3.55,1) -- (-0.5,1) -- (-1,0.5) -- (1.25,0.5);
\end{tikzpicture}
\caption{Here the parameters are $\mu=4$ and $\epsilon=0.125$, as above it was 
used Option 0 to produce it. This $2$-periodic attractor appears after the 
explosion of the $2$-periodic fat attractor in 
Figure \ref{mu4small}. As for the periodic fat attractors with 
$\epsilon=0.1115$ or $\epsilon=0.11999$, in this case the segment is also an 
attractor, its basin of attraction is like indicated at the bottom of Figure 
\ref{mu4small}.}
\label{mu4small3}
\end{figure}
This new $2$-periodic attracting set evolves towards two 
circles that we believe come from a Hopf bifurcation of an 
attracting $2$-periodic 
orbit. Pictures in Figure \ref{mu4small4} show a numerical evidence of this 
evolution.
The attracting $2$-periodic orbits in the analytic continuation of the  
$2$-periodic orbit producing that possible Hopf bifurcation are
located near $L_1$ and $L_2$, and lose 
attraction power when $\epsilon$ increases; we believe that this 
$2$-periodic orbit disappears 
and even it gives rise to attracting sets that we do not know how to identify 
or classify.

\begin{figure}[h!]
\centering
\begin{tikzpicture}
\pgftext{\includegraphics[scale=1]{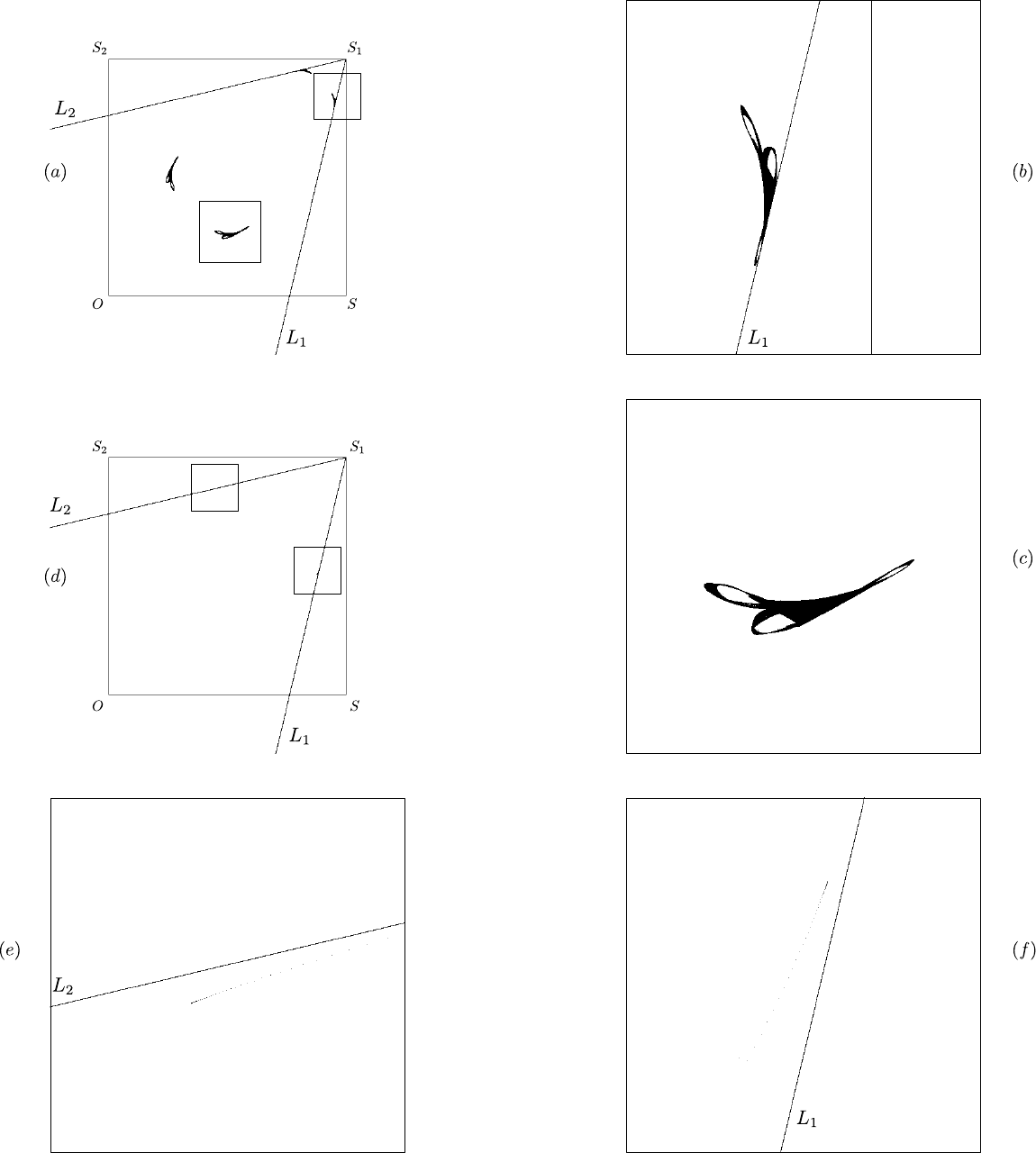}}
\draw[black!50, very thin,->] (-1.67,5.4) -- (1.22,5.4);
\draw[black!50, very thin,->] (-2.8,3.8) -- (0,3.8) -- (0,1) -- (1.2,1);
\draw[black!50, very thin,->] (-2,0.1) -- (0,0.1) -- (0,-3.4) -- (1.2,-3.4);
\draw[black!50, very thin,->] (-3.4,.68) -- (-3.4,-2.4);
\end{tikzpicture}
\caption{All the pictures in this figure were obtained applying Option 0 with 
$\epsilon=0.1915$ and $\mu=4$. This computational simulation was done with a 
lot of points, in both cases we select initial conditions whose forward orbit 
quickly accumulates in the attracting set. 
The results suggest that there are at least two attracting 
sets: an $2$-periodic orbit (located very near of $L_1$ and $L_2$) and 
apparently an $4$-periodic fat attractor. The first one is in picture $(d)$, 
indeed on the computer screen it is not clear to see it, but when magnifying 
appropriate sectors that orbit is detected; this is showed in pictures $(e)$ 
and 
$(f)$. The fat attractor is in picture $(a)$; while $(b)$ and $(c)$ are 
magnifications of the boxes indicated there.}
\label{mu4small5}
\end{figure}

\smallskip

When $\epsilon$ is between $0.14$ and $0.19$ it is relativity easy to detect 
(by means of computational simulations) an $2$-periodic orbit as the attractor 
of $F_{4,\epsilon}$ 
outside the segment $OS_1$, as above that segment is an attracting set whose 
basin is like the picture at the bottom in Figure \ref{mu4small}. 
For each $\epsilon$ with $0.19921\leq\epsilon\leq 0.25$ is observed a fat 
attractor containing the segment $OS_1$, it is unique
for the corresponding value of $\epsilon$ and it
looks like the fat attractor in Figure \ref{mu4small} with $\epsilon=0.11$, but 
with smaller area. Indeed, it is also numerically observed that when 
$\epsilon$ grows that area decreases to zero. This fact forces that the 
segment $OS_1$ is just the attracting set of $F_{4,\epsilon}$ inside 
$int\,\Gamma$ when $\epsilon$ varies from the 
value nulling that area until $\epsilon=0.5$.

\smallskip

When $\epsilon$ transits the interval $(0.19, 0.19925)$ it is really difficult 
to say globally something interesting about the attracting sets for 
$F_{4,\epsilon}$. So we only show a few pictures in Figure \ref{mu4small5} that 
are somehow a numerical confirmation of the richness and complexity of the 
dynamics of that self-map. 

\smallskip
\noindent
{\bf b) Case 2.} $\max\{4,\mu_0(\epsilon)\}<\mu\leq \mu_1(\epsilon)$\\
In this case it is clear that the basin of attraction of $\infty$ always has 
components inside $int\,\Gamma$: the triangular region $T$ described in Remark 
\ref{rk1} is always present. 
When the area of that 
components grows it is even more difficult to detect attracting sets different 
from $\infty$; indeed that area increases when $\mu$ grows, which happens if 
$\epsilon\nearrow 0.5$. Thus, we have done some computacional 
simulations with values of $\mu$ near $4$ and $\epsilon$ near those used to 
detect attracting sets when $\mu=4$; therefore it is natural that the 
attracting sets detected are, in some way, related with the dynamical 
configuration when $\mu=4$. Take for example $\epsilon=0.14$; if one increases 
the value of $\mu$, one will observe the (partial) evolution of the 
$2$-periodic orbit detected when $\mu=4$; see 
Figure \ref{mu4small4}. It can be numerically observed that this 
$2$-periodic orbit also transits by a Hopf bifurcation when $\mu$ grows from 
$4$; our computacional simulations report the following results:

\smallskip
\noindent
a) When $4\leq \mu\leq 4.00411$ the $2$-periodic orbit above seems to be the 
only attracting set of $F_{\mu,\epsilon}$; it is also noticeable the variation 
of the nature of the eigenvalues associated to that periodic orbit.

\smallskip
\noindent
b) When $\mu$ runs the interval $[4.00412,4.05284694]$, instead of the previous 
orbit there is an $2$-periodic attracting set constituted by two circles; 
between $4.00411$ and $4.00412$ a Hopf bifurcation has occurred.

\begin{figure}[h!]
\centering
\begin{tikzpicture}
\pgftext{\includegraphics[scale=1]{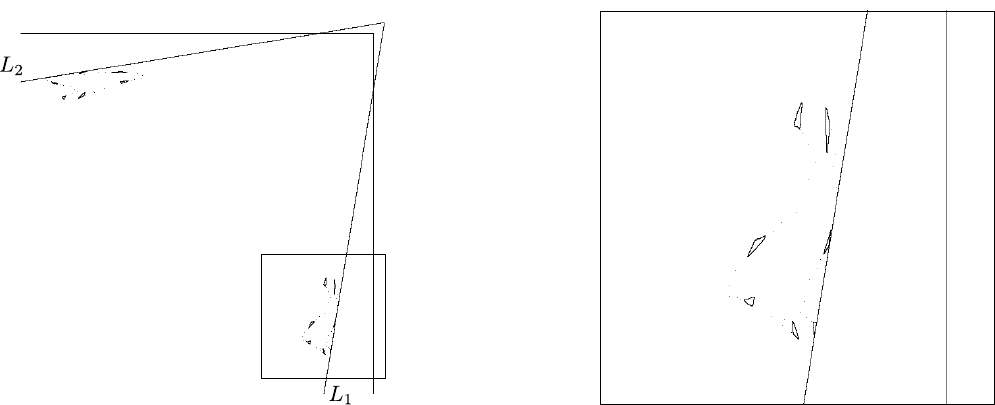}}
\draw[black!50, very thin,->] (-1.1,-1.5) -- (1,-1.5);
\end{tikzpicture}
\caption{An illustration of the periodic circles that have emerged after a 
second Hopf bifurcation when $\mu$ increases from $4$ and $\epsilon=0.14$. In 
this figure $\mu=4.085$.}
\label{caso2a}
\end{figure}

\begin{figure}[h!]
\centering
\begin{tikzpicture}
\pgftext{\includegraphics[scale=1]{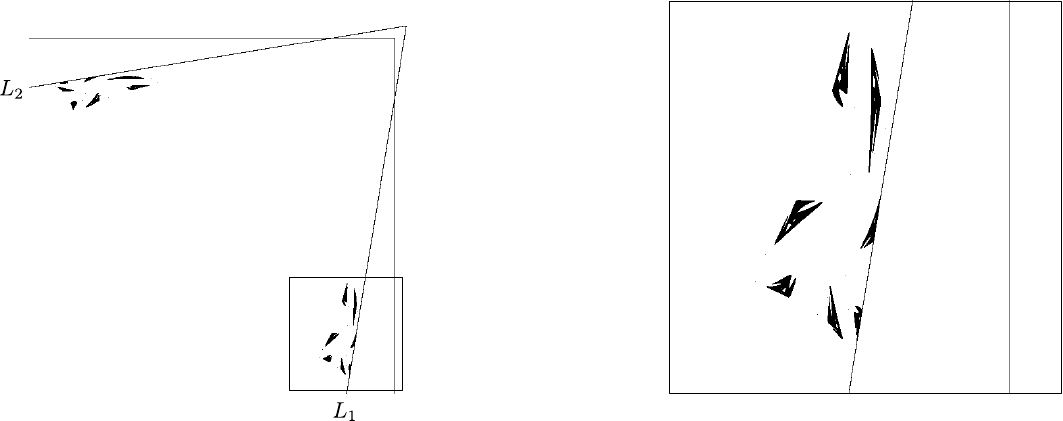}}
\draw[black!50, very thin,->] (-1.25,-1.5) -- (1.3,-1.5);
\end{tikzpicture}
\caption{This figure represents the periodic fat attractor that arises after 
the disappearance of the periodic circles that emerged from the second Hopf 
bifurcation.}
\label{caso2b}
\end{figure}

\begin{figure}[h!]
\centering
\begin{tikzpicture}
\pgftext{\includegraphics[scale=1]{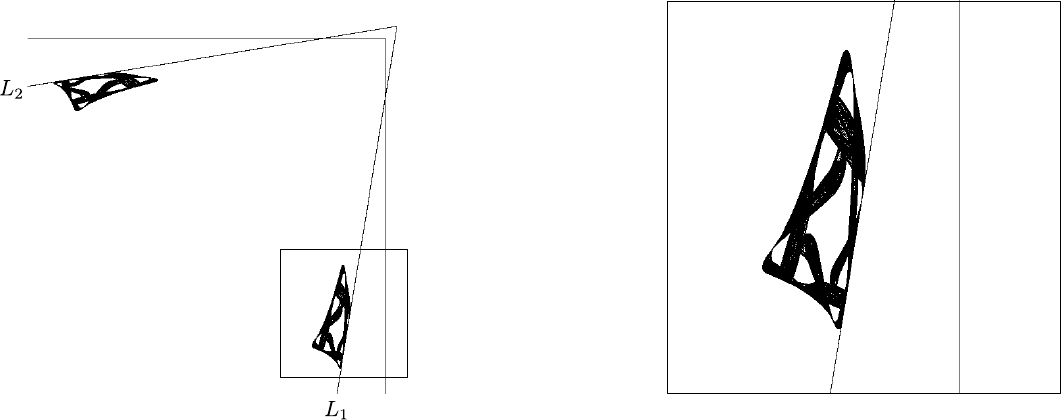}}
\draw[black!50, very thin,->] (-1.21,-1.5) -- (1.31,-1.5);
\end{tikzpicture}
\caption{This $2$-periodic fat attractor is obtained as a result of the 
collapse of the pieces (on each side of the diagonal) of the periodic fat 
attractor in the previous figure.}
\label{caso2c}
\end{figure}

\smallskip
\noindent
c) The two circles that emerged from the preceding Hopf bifurcations do not 
appear when $\mu=4.05284695$, instead there seems to be an attracting periodic 
orbit with period bigger than 2. This new periodic orbit also looks to travel 
through a Hopf bifurcation; Figure \ref{caso2a} show the 
periodic circles that arise after a new possible bifurcation, 
in this case $\mu=4.085$. These attracting periodic circles also disappear 
creating new attracting periodic orbit of even greater period. This process 
seems to repeat itself until a periodic fat attractor appear. In Figure 
\ref{caso2b} it is shown such a periodic attracting set, in this case 
$\mu=4.088$. 
In some way, of which we have no idea, the pieces
on each side of the diagonal of this 
attracting set connect to each other to form a 
piece of an $2$-periodic fat attractor; this one is shown in Figure 
\ref{caso2c} for $\mu=4.089$.

\noindent
d) The $2$-periodic fat attractors described above are easily detected 
until $\mu=4.1213$, after this value of $\mu$ we could not find, even by means 
a huge number of computational simulations, at least one attracting set. 
\begin{figure}[h!]
\centering
\includegraphics[scale=1]{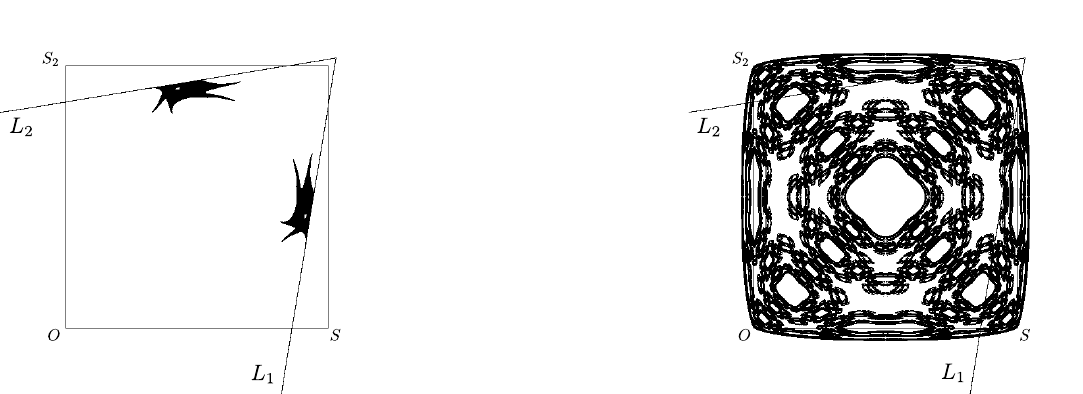}
\caption{This $2$-periodic fat attractor is obtained when $\epsilon=0.14$ and 
$\mu=4.1213$, to its right side is a graphical representation of a huge number 
of preimages of $O$, which distinguishes the components of the basin of that 
attracting set and the components of $B_\infty(F_{\mu,\epsilon})$.}
\label{caso2d}
\end{figure}
However, an obvious remark: every point in $cl(int\,\Gamma)\cap 
cl(C_{\mu,\epsilon})$ has bounded orbit, hence it has non-empty omega limit 
set, and this does not imply that there are attracting sets inside 
$int\,\Gamma$. In Figure \ref{caso2d} we show two pictures, both with 
$\epsilon=0.14$ and $\mu=4.1213$; the picture in the left side is the 
$2$-periodic fat attractor referred above, it was obtained (as always) by using 
Option $0$. The other one is the result of plotting of a large number of 
preimages of the origin, its exterior border is an approximation of 
$\Gamma$, the white sectors correspond to: the components of 
$B_\infty(F_{\mu,\epsilon})$ inside $int\,\Gamma$, and the components of the 
basin of attraction of the $2$-periodic fat attractor.

\subsubsection{\bf Invariant curves beyond 
$\boldsymbol{\epsilon\longmapsto \mu_1(\epsilon)}$}

As announced at the end of item 2 in Remark \ref{rk1}, we have no a proof for 
the existence of the curve $\Gamma$ (as introduced in Theorem \ref{a}) when 
$\mu>\mu_1(\epsilon)$ and it is close to 
$\mu_1(\epsilon)$; 
\begin{figure}[h!]
\centering
\begin{tikzpicture}
\pgftext{\includegraphics[scale=1]{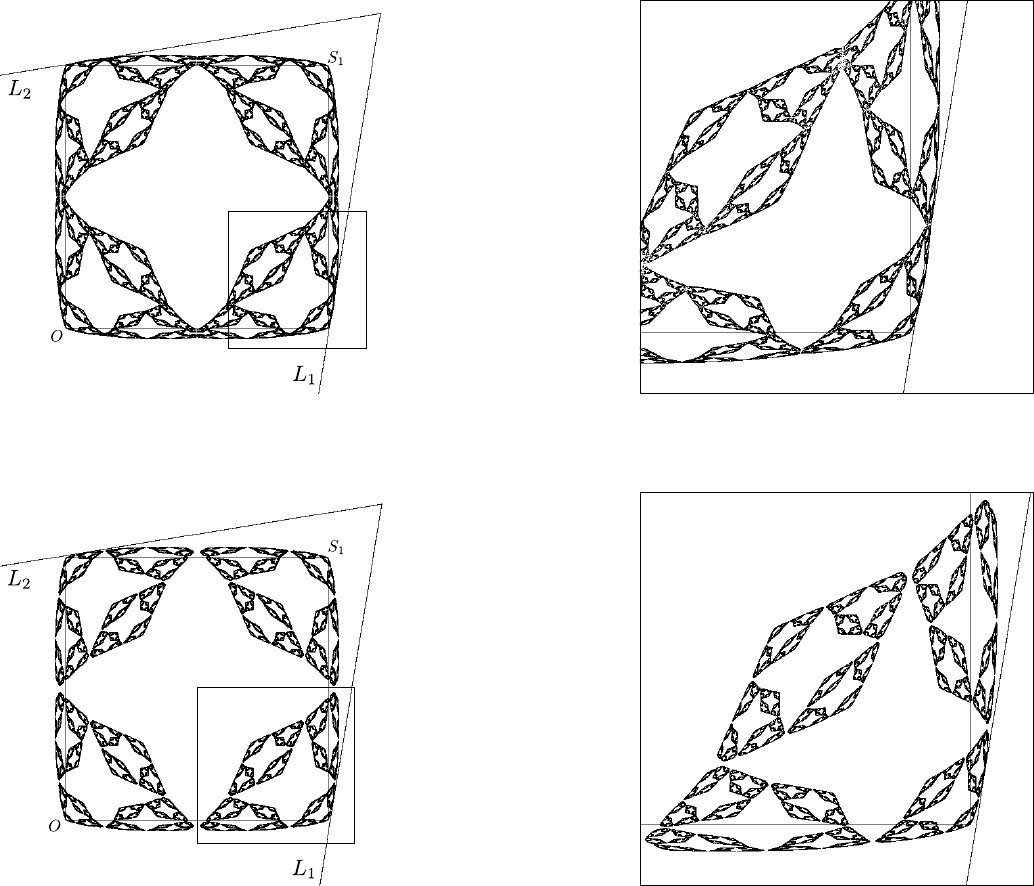}}
\draw[black!50, very thin,->] (-1.6,-3.5) -- (1.2,-3.5);
\draw[black!50, very thin,->] (-1.48,1.5) -- (1.2,1.5);
\end{tikzpicture}
\caption{For the pictures at the top we have used $\epsilon=0.14$ and 
$\mu=4.8$, they show an invariante curve $\Gamma$ beyond $\epsilon\longmapsto 
\mu_1(\epsilon)$; while for the pictures in lower part $\epsilon=0.14$ and 
$\mu=4.82$, with these values there is no curve as stated in 
Theorem \ref{a}.}
\label{beyond}
\end{figure}
however, there we 
discussed some ideas of how such a proof could be made.
Pictures at the top of Figure \ref{beyond} numerically show, for $\epsilon=0.14$ 
and $\mu=4.8$,  
the existence of a positively invariant curve $\Gamma$ (exterior 
border of picture in the left side)
such that $ext\,\Gamma$ is the immediate basin of attraction of $\infty$. For
in this computational simulation the point $S$ is in the interior 
of the cone $C_{\mu,\epsilon}$; that is $\mu>\mu_1(\epsilon)$, see the 
magnification of the box in the picture located in the upper left corner. 
At this time we consider appropriate to establish the following conjecture:

\smallskip
\noindent
{\bf Conjecture:}
{\em If $\epsilon\in (0,1/2)$ and $\mu>1$ satisfy $\mu>\mu_1(\epsilon)$ and 
there exists a curve $\Gamma$ as in Theorem \ref{a}, then 
$B^c_\infty(F_{\mu,\epsilon})$ is the closure of 
$\bigcup_{n\geq 0}F_{\mu,\epsilon}^{-n}(\Gamma)$.}

\smallskip

In the small strength case, the computational simulations reveal that 
there is always an invariant curve $\Gamma$ as in Theorem \ref{a}
whenever
$\mu>\mu_1(\epsilon)$ and 
close enough to $\mu_1(\epsilon)$. We would like to recall that 
this curve disappears when the value of 
$\mu$ grows, this certainly occurs when from a certain 
preimage of the circle $C$, the following ones are inside $C_{\mu,\epsilon}$. 
In this case the basin of attraction 
$B_\infty(F_{\mu,\epsilon})$ is again connected and its complement has 
infinitely many components. It is well known that this last set is an expanding 
Cantor set for all $\mu$ large enough, see \cite{rrv}.
\subsection{Large strength case}\label{large}
The aim of this part is analagous to that of the subsection \ref{small}.
In this setting it is well known that if $\mu$ and $\epsilon$ satisfy  
$1<\mu\leq\mu_0'(\epsilon)$, 
then the set of points with bounded orbit is the union of the segments 
$OS_1$ and $SS_2$; in addition, the forward orbit of
every $z$ in this union converges to the fixed point 
$P_\mu=(\frac{\mu-1}{\mu},\frac{\mu-1}{\mu})$. On the other hand, if 
$1<\mu\leq\mu_1(\epsilon)$, then the phenomenon of 
Cantor's set of circles discussed in the small straight case is not possible, 
this is due to the absence of a Cantor set on the diagonal. 

\subsubsection{\bf Attractors inside $\boldsymbol{int\,\Gamma}$}
As in the small strength case, there are regions in the parameter space with 
$1<\mu\leq\mu_1(\epsilon)$ where the map $F_{\mu,\epsilon}$ exhibits fat 
attractors when $\mu$ and $\epsilon$ 
are near $4$ and $0$, respectively. In order to not to be repetitive we do not 
show such examples; however, we will show some type of attracting sets that we 
did 
not exhibit before.

\smallskip

\begin{figure}[h!]
\centering
\includegraphics[scale=1]{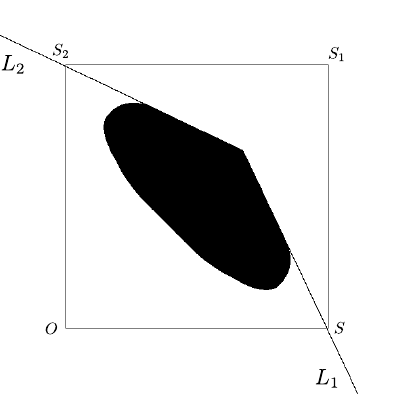}
\hspace{1cm}
\includegraphics[scale=1]{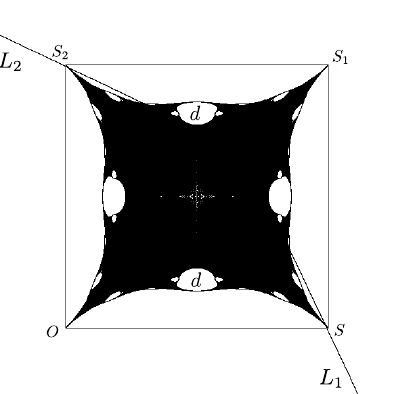}
\caption{Fat attractor (left picture) and its basin (right picture) when 
$\epsilon=-0.9$ and 
$\mu=2.71$.}
\label{fat3}
\end{figure}

\begin{figure}[h!]
\centering
\includegraphics[scale=1]{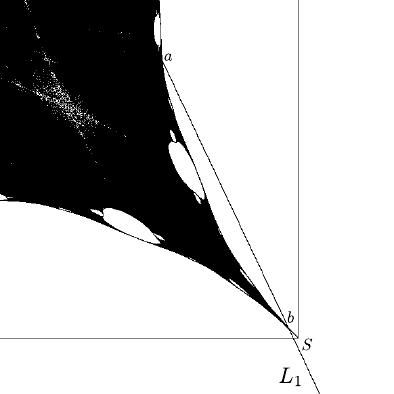}
\caption{The sector whose boundary is given by the segment $ab$ in $L_1$ and 
the arc in $\Gamma$ between the points $a$ and $b$ produces some components 
of $B_\infty(F_{\mu,\epsilon})$ intersecting the critical line $\ell_1$.}
\label{fat4}
\end{figure}

Take $\epsilon=-0.9$, so $\mu_1(\epsilon)\sim 2.714$ and 
$\mu_2(\epsilon)\sim 2.357$. For $\mu=2.71$ 
the map $F_{\mu,\epsilon}$ has the fat attractor plotted in Figure 
\ref{fat3} (left picture). 
This attractor has been generated applying Option 0 to the point 
$(0.67,0.59)$, just with this point and Option 2 one obtains an illustration of 
the basin of that attractor; see picture on the right side in Figure 
\ref{fat3}. Observe that the exterior border of that picture is an 
approximation of the Jordan curve $\Gamma$.
We take advantage of it to reiterate the existence of components 
(white 
regions) of 
$B_\infty(F_{\mu,\epsilon})$ inside $int\,\Gamma$, this was 
discussed in item 2 
of Remark \ref{rk2}. Recall that there exists only one way to 
produce these components in the large strength case. In 
Figure \ref{fat4} we show a magnification of a 
neighborhood of the point $S$ where it can be seen an arc of $\Gamma$ through 
which are produced the two components of $B_\infty(F_{\mu,\epsilon})$ indicated 
by $d$ in the picture in the right side of Figure \ref{fat3}.

\begin{figure}[h!]
\centering
\includegraphics[scale=1]{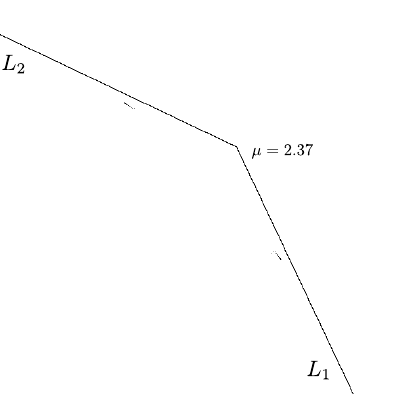}
\includegraphics[scale=1]{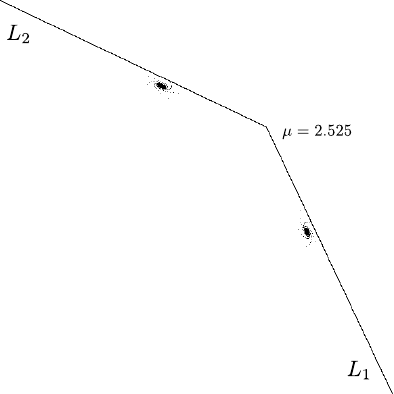}
\includegraphics[scale=1]{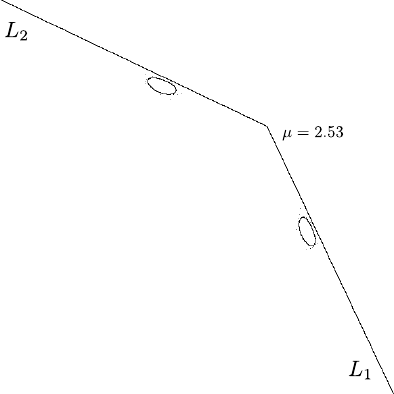}

\vspace*{1cm}
\includegraphics[scale=1]{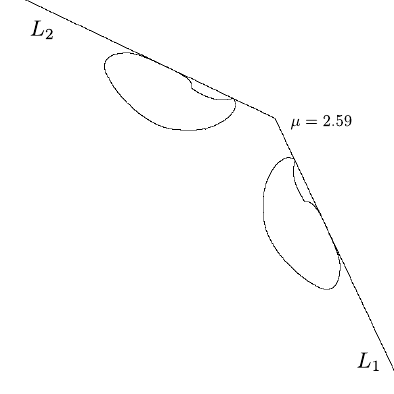}
\includegraphics[scale=1]{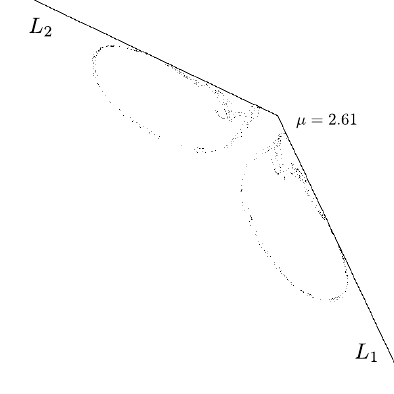}
\includegraphics[scale=1]{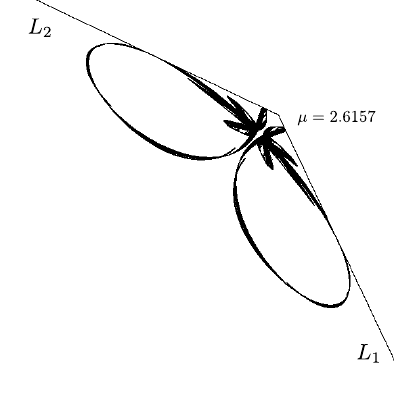}

\vspace*{1cm}
\includegraphics[scale=1]{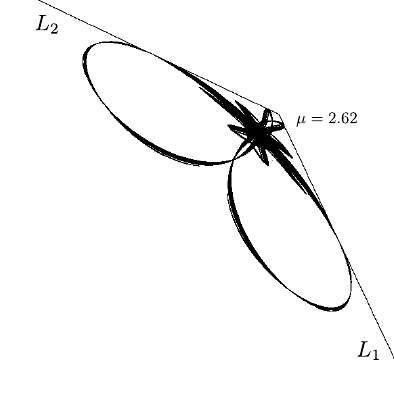}
\includegraphics[scale=1]{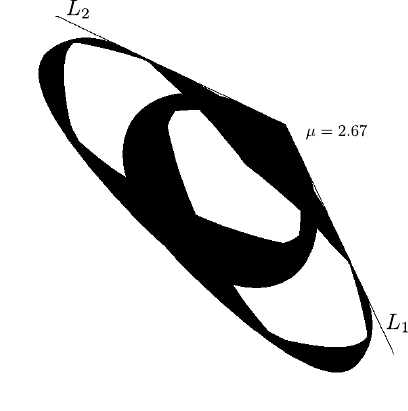}
\includegraphics[scale=1]{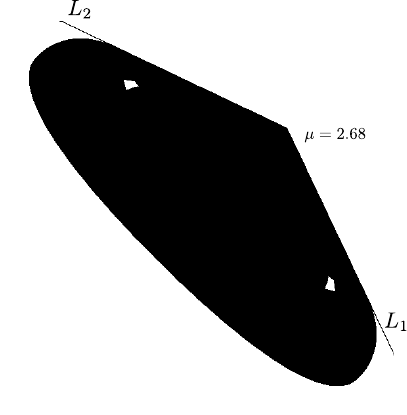}
\caption{A simplistic and inaccurate graphical illustration of the evolution 
towards the fat attractor in Figure \ref{fat3}.}
\label{fat5}
\end{figure}

Try to explain, even to understand, how the fat attractor in Figure \ref{fat3} 
evolved when $\mu$ varies and $\epsilon$ stays equal is a very hard task. In 
the following pictures we present some numerical experiments, with 
$\epsilon=-0.9$, which partially show the possible attractors of 
$F_{\mu,\epsilon}$ inside $int\,\Gamma$. In that picture sequence the 
parameter $\mu$ takes some values in the interval 
$(\mu_2(\epsilon),\mu_1(\epsilon))$ in an increasing way. We recall that if 
$\mu\in (\mu'_0(\epsilon),\mu_2(\epsilon)]$, the fixed point $P_\mu$ is the 
only attractor of $F_{\mu,\epsilon}$ inside $int\,\Gamma$.
It is important to say that these numerical experiments have been repeated for a 
large number of values of $\epsilon<0$ ​​and similar figures or attractors have 
been observed to those shown in Figure \ref{fat5}.
In the picture corresponding to $\mu=2.37$ it is observed as attractor an 
$2$-periodic orbit. From the numerical point of view, this dynamical 
configuration 
remains for a certain range of values of $\mu$ in which the nature of the 
eigenvalues associated to that periodic orbit evolves towards a Hopf 
bifurcations, which should occur for some  
$2.525<\mu<2.53$. The picture associated to $\mu=2.525$ induces to think that 
the eigenvalues of the $2$-periodic orbit are non-real with negative real part; 
the picture corresponding to
$\mu=2.53$ shows the $2$-periodic attracting circles generated 
by that Hopf bifurcation. These two circles go through a degenerative stage 
losing differentiability, approaching to the boundary of $C_{\mu,\epsilon}$ 
until they disperse in a cloud of points (perharps a periodic orbit of very 
large period) to later constitute a fat periodic attractor; this part of the 
evolution is briefly illustrated in the pictures corresponding to 
$\mu=2.59$, $\mu=2.61$ and $\mu=2.6157$. Finally, the pictures associated to 
$\mu=2.62$, $\mu=2.67$ and $\mu=2.68$ show the collapse of that fat periodic 
attractor in a single one, gaining area until reaching the fat 
attractor in Figure \ref{fat3}. 

\subsubsection{\bf Invariant curves beyond 
$\boldsymbol{\epsilon\longmapsto \mu_1(\epsilon)}$}
\begin{figure}[h!]
\centering
\includegraphics[scale=1]{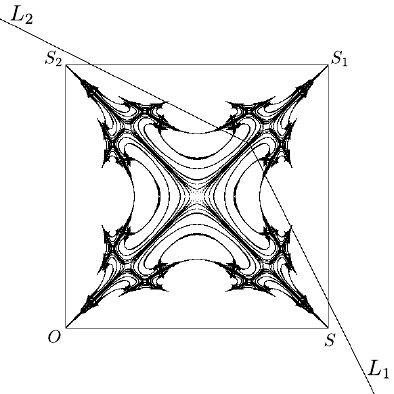}\hspace{1cm}
\includegraphics[scale=1]{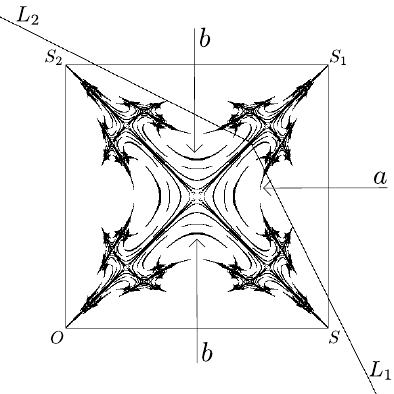}
\caption{Picture on the left side was generated with $\epsilon=-1$ and 
$\mu=2.8$, while the one on the right side was obtained with $\epsilon=-1$ and 
$\mu=2.82$.}
\label{nodif}
\end{figure}

In the large strength case
when $\mu_1(\epsilon)<\mu<4$ the preimage of 
$\partial Q$ continues to be a closed curve joining the points in 
$F^{-1}_{\mu,\epsilon}(O)$; so it is still possible to discuss 
the existence of a forward invariant Jordan curve containing 
$F^{-1}_{\mu,\epsilon}(O)$; see item 3 of Remark \ref{rk2}. 
However, one can not hope neither differentiability of the arcs 
$\Gamma_b,\Gamma_\ell,\Gamma_r$ and $\Gamma_t$ of $\Gamma$, nor 
$ext\,\Gamma\subset B_\infty(F_{\mu,\epsilon})$; pictures in Figure \ref{nodif}
are numerical samples of these possible features.
Since the points $S$ and $S_2$ have preimages, we have used Option 1 (applied 
to the origin) to produce that pictures; notice that their exterior borders 
look likes non-smooth closed curves. In the picture on the right 
side it can be seen that the preimage of the set 
indicated with the letter $a$ contains points with bounded orbirts, which 
belong to the set indicated with $b$; recall what was discussed at the end of 
item 3 in 
Remark \ref{rk2}.

\smallskip

\begin{figure}[h!]
\centering
\includegraphics[scale=1]{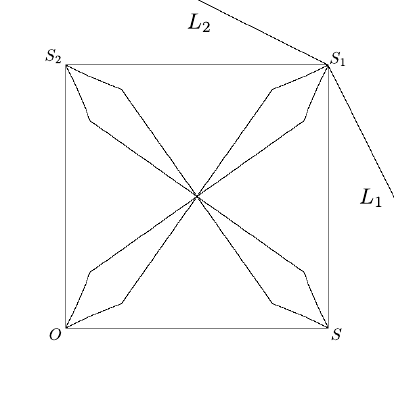}\hspace{1cm}
\includegraphics[scale=1]{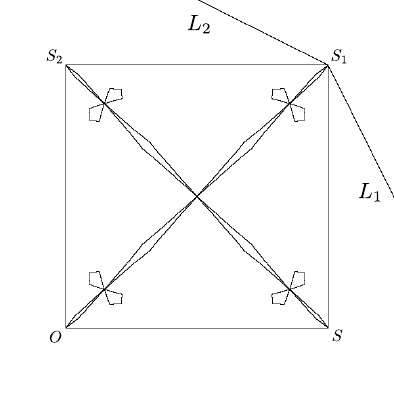}
\caption{These pictures represent the two first preimages of $\partial Q$ when 
$\epsilon=-1$ They are similar for all negative value of $\epsilon$. }
\label{nodif2}
\end{figure}

Now we discuss the special case $\mu=4$. Observe that
the square $Q$ is contained in the interior of $C_{\mu,\epsilon}$ except 
by the point $S_1$, which is just the vertex of that cone; so, if there 
exists an invariant Jordan curve $\Gamma$ as described in Theorem \ref{b}, then 
there are no components of $B_\infty(F_{\mu,\epsilon})$ inside $int\,\Gamma$. 
Also 
observe that the preimage of $Q$ is union of four compact set $Q_i$ 
($i=0,1,2,3$) such that $Q_i\cap Q_j=\{(1/2,1/2)\}$ for all $i\neq j$. In 
particular $F^{-1}_{\mu,\epsilon}(\partial Q)$ is the union of four closed 
curves, one for each point in $F^{-1}_{\mu,\epsilon}(O)$ and having $(1/2,1/2)$ 
as the only common point; pictures in Figure \ref{nodif2} show  
$F^{-1}_{\mu,\epsilon}(\partial Q)$ and $F^{-2}_{\mu,\epsilon}(\partial Q)$, 
which were obtained with Option 4 and $\epsilon=-1$.
It is clear that one can repeat the procedure discussed in item 3 of Remark 
\ref{rk2} to construct a curve $\Gamma$. Another way to achieve the same goal 
is to resort to the tools that the Iterated Function Systems theory 
(\cite{barnsley}, \cite{hutchison}) provides. The corresponding iterated 
function system acts on  
the compact metric space $Q$ and it is given by the four inverse branches $H_i$ 
($i=0,1,2,3$) 
of $F_{\mu,\epsilon}$, which are defined by means of 
\eqref{pre}. The aim of this procedure is to show that the iteration of the 
Hutchinson operator related to that iterated function system converges to a 
unique attractor: the curve $\Gamma$. This property is satisfied if, for 
example, each self-map 
$H_i:Q\to Q$ is 
such that one of its iterated is a contraction, a task that in this case looks a 
bit difficult. In any case, in the following figure we show an approximation of 
$B^c_\infty(F_{\mu,\epsilon})$ when $\epsilon=-1$; it has been produced 
plotting a huge number of preimages of the origin. We highlight that this 
self-similar fractal set is analogous to that generated by the same numerical 
experiment with any value of $\epsilon<0$.
\begin{figure}[h!]
\centering
\includegraphics[scale=1]{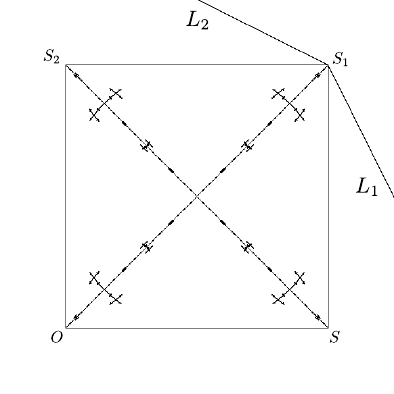}\hspace{1cm}
\includegraphics[scale=1]{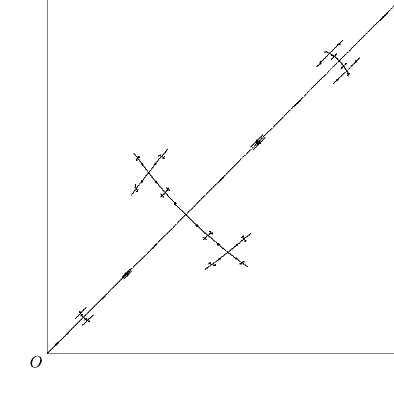}
\caption{Illustration of the self-similar fractal set representing the set of 
points with bounded orbits when $\mu=4$. Picture on the right side is a zoom 
near the origin.}
\label{nodif3}
\end{figure}

\medskip
\noindent
{\bf Acknowledgments.} This research was suppoted by Grant 002-CT-2015 from 
the Consejo de Desarrollo Científico, Humanístico y Tecnológico (CDCHT) of the 
Universidad Centroccidental Lisandro Alvarado.

\end{document}